\magnification=1200
\hoffset=.25truein \advance\hsize by -0.5truein
\def\leaderfill{\leaders\hbox to 0.5em{\hss.\hss}\hfill}
\def\itc#1{{\it #1\/}}
\def\X{{\rm X}}
\font\titlefont=cmbx12 scaled \magstep1

{
\nopagenumbers
\def\leaderfill{\leaders\hbox to 0.6em{\hss.\hss}\hfill}
\parindent=0pt
\baselineskip=16.8pt
\font\titlefont=cmbx12 scaled \magstep2

\newskip\mdskipamount
\mdskipamount=7.2pt plus2.4pt minus2.4pt
\def\mdskip{\vskip\mdskipamount}
\newskip\bskipamount
\bskipamount=14.4pt plus4.8pt minus4.8pt
\def\bskip{\vskip\bskipamount}

\font\rm=cmr12
\rm

\centerline{\titlefont Symmetric Function Generalizations}
\mdskip
\centerline{\titlefont of Graph Polynomials}
\mdskip
\centerline{\rm by}
\centerline{\rm Timothy Yi-Chung Chow}
\mdskip
\centerline{\rm A.B., Princeton University, New Jersey (1991)}
\mdskip
\centerline{\rm Submitted to the Department of Mathematics}
\centerline{\rm in partial fufillment of the requirements for the degree of}
\mdskip
\centerline{\rm Doctor of Philosophy}
\mdskip
\centerline{\rm at the}
\mdskip
\centerline{\rm MASSACHUSETTS INSTITUTE OF TECHNOLOGY}
\mdskip
\centerline{\rm June 1995}
\mdskip
\centerline{\rm \copyright 1995, Timothy Yi-Chung Chow
}
\mdskip
\centerline{\rm The author hereby grants to MIT permission to reproduce and}
\centerline{\rm to distribute copies of this thesis document in whole or in
             part.}
\bskip\bskip
\line{Signature of Author\leaderfill}
\rightline{Department of Mathematics}
\rightline{April 26, 1995}
\bskip\bskip
\line{Certified by\leaderfill}
\rightline{Richard P. Stanley}
\rightline{Professor of Applied Mathematics}
\rightline{Thesis Supervisor}
\bskip
\line{Accepted by\leaderfill}
\rightline{David Vogan}
\rightline{Chairman, Departmental Graduate Committee}
\rightline{Department of Mathematics}
\vfill\eject
}

{\nopagenumbers
\hbox{ }
\vfill\eject}

\centerline{\titlefont Symmetric Function Generalizations}
\medskip
\centerline{\titlefont of Graph Polynomials}
\medskip
\centerline{by}
\centerline{Timothy Yi-Chung Chow}
\bigskip
\centerline{Submitted to the Department of Mathematics}
\centerline{on April 26, 1995, in partial fulfillment of the}
\centerline{the requirements for the degree of}
\centerline{Doctor of Philosophy}
\bigskip\bigskip
\leftline{\titlefont Abstract}
\medskip
\noindent
Motivated by certain conjectures regarding
immanants of Jacobi-Trudi matrices,
Stanley has recently defined and
studied a symmetric function generalization~$\X_G$
of the chromatic polynomial of
a graph~$G$.  Independently, Chung
and Graham have defined and
studied a directed graph invariant
called the cover polynomial.  The
cover polynomial is closely related
to the chromatic polynomial and
to the rook polynomial, so
it is natural to ask
if one can mimic Stanley's
construction and generalize the cover
polynomial to a symmetric function.
The answer is yes, and
the bulk of this thesis
is devoted to the study
of this generalization, which we
call the path-cycle symmetric
function.  We obtain analogues of
some of Stanley's theorems about~$\X_G$
and we generalize some of
the theory of the cover
polynomial and the rook polynomial.
In addition, we are led
to define a symmetric function
basis that seems to be
a ``natural'' generalization of the
polynomial basis $$\biggl\{{x+k\choose d}\biggr\}_{k=0,1,\ldots,d}$$
and we prove a combinatorial
reciprocity theorem that gives an
affirmative answer to Chung and
Graham's question of whether the
cover polynomial of a digraph
determines the cover polynomial of
its complement.  The reciprocity
theorem also ties together several
scattered results in the literature
that previously seemed unrelated.

In the remainder of the
thesis, we prove some miscellaneous
results about Stanley's function~$\X_G$
and we also sketch briefly
in the introduction a possible
approach to generalizing other graph
polynomials (and a few other
combinatorial polynomials) to symmetric functions.

\bigskip
\leftline{Thesis supervisor: Richard P. Stanley}
\leftline{Title: Professor of Applied Mathematics}
\vfill\eject
{
\nopagenumbers
\hbox{ }
\vfill\eject
\font\dedfont=pzcmi at 18pt
\vglue 2.5in
\centerline{\dedfont To my parents}
\vfill\eject
\hbox{ }
\vfill\eject
}

\parskip=10pt plus4pt minus4pt

\font\titlefont=pplb at 20pt
\vbox to 1.5truein{\vfill}\leftline{\titlefont Acknowledgments}\vbox
        to 1.5truein{\vfill}

I am deeply grateful to my thesis advisor, Richard Stanley.
Without him, not only would this thesis not exist, but I would be
entirely ignorant of this beautiful area of combinatorics (and many
other areas besides) to this day.  Moreover, while many graduate
students have a long stream of advisor horror stories,
and regard their graduate school experience as an
agony in eight fits, I have no complaints at all about Richard.
My experience here at M.I.T. has been remarkably pleasant, with
the main source of stress being my frustration at my own ignorance
and laziness.

I also wish to thank Ira Gessel and Gian-Carlo Rota, who have been
consistently encouraging and unreasonably patient with my frequent
silly questions and comments.  I am also indebted to John Stembridge,
whose SF and posets packages for Maple were of immeasurable help to
me while I was doing my research.

Seeing so many other students struggling to pay for their education,
I am reminded how fortunate I am to enjoy
the generous financial support of the M.I.T. mathematics department
(through teaching assistantships) and the National Science Foundation 
(through graduate fellowship stipends), who have made it possible
for me to pursue my studies without having to worry about money.

To my family and my friends from Princeton, Duluth, Chinese
Christian Fellowship, Graduate Christian Fellowship, and the
math department---to all who have loved me, lived with me,
taught me, hung out with me, laughed with me, tolerated me,
encouraged me, and prayed for me---I apologize for not
expressing my appreciation to you more often.  Without you,
the past four years would have been bleak indeed, and I would
have been lucky to survive at all, let alone write a thesis
and graduate from M.I.T.

Last of all, and most of all, I thank God for giving me life,
filling it with joy, and creating this thing of beauty that we
call mathematics.
\vfill\eject

{\nopagenumbers
\hbox{ }
\vfill\eject}

{
\baselineskip=16pt
\parindent=0pt
\def\ind{\kern 1truein}

\vbox to 1.5truein{\vfill}\leftline{\titlefont Contents}\vbox
        to 1.5truein{\vfill}

\line{{\bf Acknowledgments}\leaderfill 7}
\line{{\bf Introduction}\leaderfill 11}
\line{{\bf Chapter 1: Preliminaries}\leaderfill 19}
\line{\ind 1. Graphs, digraphs, and symmetric functions\leaderfill 19}
\line{\ind 2. $\chi_G$, $\X_G$, $C(D)$, and $\Xi_D$\leaderfill 21}
\line{\ind 3. Connection with rook theory\leaderfill 23}
\line{{\bf Chapter 2: The Path-Cycle Symmetric Function}\leaderfill 24}
\line{\ind 1. Basic facts\leaderfill 24}
\line{\ind 2. Reciprocity\leaderfill 30}
\line{\ind 3. Rook theory\leaderfill 38}
\line{\ind 4. The Poset Chain Conjecture\leaderfill 49}
\line{{\bf Chapter 3: The Chromatic Symmetric Function}\leaderfill 53}
\line{\ind 1. $G$-ascents\leaderfill 53}
\line{\ind 2. The Poset Chain Conjecture revisited\leaderfill 55}
\line{\ind 3. Reconstruction\leaderfill 57}
\line{\ind 4. Superfication\leaderfill 59}
\line{\ind 5. $\X_G(t)$\leaderfill 64}
\line{\ind 6. Counterexamples\leaderfill 66}
\line{{\bf Bibliography}\leaderfill 67}
\vfill\eject
}

{\nopagenumbers
\hbox{ }
\vfill\eject}

\catcode`\@=11
 
\def\hexnumber@#1{\ifnum#1<10 \number#1\else
 \ifnum#1=10 A\else\ifnum#1=11 B\else\ifnum#1=12 C\else
 \ifnum#1=13 D\else\ifnum#1=14 E\else\ifnum#1=15 F\fi\fi\fi\fi\fi\fi\fi}
 
\font\tenmsa=msam10
\font\sevenmsa=msam7
\font\fivemsa=msam5
\newfam\msafam
\textfont\msafam=\tenmsa  \scriptfont\msafam=\sevenmsa
  \scriptscriptfont\msafam=\fivemsa
\def\msa@{\hexnumber@\msafam}

\mathchardef\surject="3\msa@10
\mathchardef\therefore="3\msa@29
\mathchardef\because="3\msa@2A
\mathchardef\normalsup="3\msa@42
\mathchardef\normalsub="3\msa@43
\mathchardef\commutes="\msa@09
\mathchardef\curlyle="3\msa@2D
\mathchardef\curlyge="3\msa@25
\mathchardef\revreac="3\msa@0A
\mathchardef\varle="3\msa@36
\mathchardef\varge="3\msa@3E
\mathchardef\yen="0\msa@55
\mathchardef\varpropto="3\msa@5F
\mathchardef\regtm="0\msa@72
\mathchardef\plusdot="2\msa@75
 
\font\tenmsb=msbm10
\font\sevenmsb=msbm7
\font\fivemsb=msbm5
\newfam\msbfam
\textfont\msbfam=\tenmsb  \scriptfont\msbfam=\sevenmsb
  \scriptscriptfont\msbfam=\fivemsb
\def\msb@{\hexnumber@\msbfam}

\mathchardef\nodiv="3\msb@2D
\mathchardef\coveredby="3\msb@6C
\mathchardef\covers="3\msb@6D
\mathchardef\propsub="3\msb@28
\mathchardef\propsup="3\msb@29
\mathchardef\varemptyset="0\msb@3F

%
%

\font\teneuf=eufm10
\font\seveneuf=eufm7
\font\fiveeuf=eufm5
\newfam\euffam
\textfont\euffam=\teneuf \scriptfont\euffam=\seveneuf
  \scriptscriptfont\euffam=\fiveeuf
\def\euf@{\hexnumber@\euffam}
\def\Frak{\relax\ifmmode\let\next\Frak@\else
 \def\next{\errmessage{Use \string\Frak\space only in math mode}}\fi\next}
\def\Frak@#1{{\Frak@@{#1}}}
\def\Frak@@#1{\fam\euffam#1}

\font\tenmsb=msbm10
\font\sevenmsb=msbm7
\font\fivemsb=msbm5
\newfam\msbfam
\textfont\msbfam=\tenmsb  \scriptfont\msbfam=\sevenmsb
  \scriptscriptfont\msbfam=\fivemsb
\def\msb@{\hexnumber@\msbfam}
\def\Bbb{\relax\ifmmode\let\next\Bbb@\else
 \def\next{\errmessage{Use \string\Bbb\space only in math mode}}\fi\next}
\def\Bbb@#1{{\Bbb@@{#1}}}
\def\Bbb@@#1{\fam\msbfam#1}

%

\font\tenscr=rsfs10 
\font\sevenscr=rsfs7 
\font\fivescr=rsfs5 
\skewchar\tenscr='177 \skewchar\sevenscr='177 \skewchar\fivescr='177
\newfam\scrfam
\textfont\scrfam=\tenscr \scriptfont\scrfam=\sevenscr
  \scriptscriptfont\scrfam=\fivescr
\def\Scr{\relax\ifmmode\let\next\Scr@\else
 \def\next{\errmessage{Use \string\Scr\space only in math mode}}\fi\next}
\def\Scr@#1{{\Scr@@{#1}}}
\def\Scr@@#1{\fam\scrfam#1}

\catcode`\@=12


\def\defeq{\buildrel \rm def \over =}
\parskip=12pt plus4pt minus4pt

\vbox to 1.5truein{\vfill}\leftline{\titlefont Introduction}\vbox
        to 1.5truein{\vfill}

Mathematicians are notorious for sterilizing their technical papers so
that there remains no taint of motivation.  All traces of how the proofs
were originally discovered are completely obliterated---or, at the very
least, carefully disguised---so that maximum brevity and elegance are
achieved.  While this approach is ideal for the reader
who wishes to use the paper solely as a reference (or as an object for
aesthetic contemplation), the student who really wants to understand the
paper is done a disservice.  On the other hand, an informal presentation
that gives a great deal of intuitive motivation and works through many
concrete examples may be ideal pedagogically, but it is frustrating for
the researchers who wish to extract only what they need for their
own~work.

I have adopted a compromise here.  The body of this thesis is written
in the customary definition-theorem-proof style, with some
motivational comments scattered haphazardly.
However, this introduction is written in a
more informal style, and I try to show how my work
fits into the larger scheme of things.
I provide precise definitions here only when I feel that
they are critical to the exposition; for full definitions the reader
should consult Chapter~1 and the references cited therein.

The principal objects of study in this thesis may not seem natural or
intrinsically interesting a priori, so let me begin by describing two
previous lines of research that \itc{are} interesting and natural.
We will then see how the present thesis arose out of the unexpected
convergence of those two lines of study.

The first of these pre-existing lines of research is the study
of \itc{immanants} of matrices.  To understand what these are, recall
that the \itc{determinant} of a matrix may be thought of as
a sum of the form
$$\sum_{\sigma\in S_n} (-1)^\sigma \prod_{i=1}^n x_{i,\sigma(i)},$$
where $S_n$ denotes the symmetric group on $n$~letters.
Similarly, the permanent of a matrix is a sum of the form
$$\sum_{\sigma\in S_n} \prod_{i=1}^n x_{i,\sigma(i)}.$$
Notice that the only difference between these two expressions lies in
the coefficients: in the case of the determinant, the coefficients are
given by the sign character of~$S_n$, and in the case of the
permanent, the coefficients are given by the trivial character
of~$S_n$.  This observation suggests that perhaps expressions of
the form
$$\sum_{\sigma\in S_n} \chi^\lambda(\sigma)
   \prod_{i=1}^n x_{i,\sigma(i)}$$
(where $\chi^\lambda$ is an irreducible character of~$S_n$), called
\itc{immanants,} might be interesting objects of study.  It turns out,
not entirely unexpectedly, that the theory of immanants is not as rich
as the theory of determinants (see [Gre] for a brief survey of work done
on immanants).  However, in the special case of \itc{Jacobi-Trudi
matrices,} which are matrices that arise in the theory of symmetric
functions, there is a rich array of combinatorial conjectures related
to immanants.  While it is tempting to give a summary of these beautiful
conjectures, we shall restrict ourselves to describing just one of them,
due to Stanley and Stembridge~[S-S], since it is the only one directly
relevant to our present purposes.  The interested reader is referred
to~[Gre] and~[S-S] for the full story.

To state the Stanley-Stembridge conjecture (also known as the \itc{Poset
Chain Conjecture}), we must first define a certain invariant~$\X_P$ of a
poset~$P$.  A \itc{coloring} of a poset~$P$ is a map $\kappa$ that
sends each element of~$P$ to a positive integer (or ``color'') and
whose fibers are totally ordered subsets (or ``chains'') of~$P$.
Then $\X_P$ is defined~by
$$\X_P \defeq \sum_\kappa \prod_{p\in P} x_{\kappa(p)},$$
where the sum is over all colorings of~$P$.
Since the colors may be permuted indiscriminately, $\X_P$ is a symmetric
function, i.e., it is invariant under any
permutation of the indeterminates~$x_i$.
By the well-known fundamental theorem
of symmetric functions, every symmetric function has a unique expression
as a polynomial in the \itc{elementary symmetric functions~$e_n$}
defined by
$$e_n = \sum_{0<i_1<i_2<\cdots<i_n} x_{i_1} x_{i_2} \cdots x_{i_n}.$$
The Poset Chain Conjecture states that if $P$ does not contain an
induced subposet isomorphic to the disjoint union of a three-element
chain and a one-element chain, then $\X_P$ has nonnegative coefficients
when expressed as a polynomial in the elementary symmetric functions.
(For the exact connection between this conjecture and the immanants of
Jacobi-Trudi matrices, the reader is referred to~[S-S].)
The Poset Chain Conjecture is supported by considerable numerical evidence
and some partial results (notably Gasharov's theorem~[Ga1] that under the
hypotheses of the conjecture, $\X_P$ has nonnegative coefficients when 
expanded in terms of Schur functions), but remains open.

The next step in the story is a simple but important
observation of Stanley that the invariant~$\X_P$
can be extended from a poset invariant to a graph invariant.
More specifically, let $G(P)$ denote the
\itc{incomparability graph} of~$P$, i.e., the graph whose vertex set
is the set of elements of~$P$ and in which two vertices are adjacent
if and only if they are \itc{not} comparable as elements of~$P$.
Note that chains in~$P$ correspond to independent sets of~$G(P)$
(i.e., subsets of vertices with no edges between them),
and that colorings of~$P$ correspond to
colorings of~$G(P)$ in the usual sense of assignments of colors
to the vertices of~$G(P)$ such that adjacent vertices are never
assigned the same color.
Now the notion of a coloring makes sense for any
graph~$G$ and not just incomparability graphs, so we may define
$$\X_G \defeq \sum_\kappa \prod_{v} x_{\kappa(v)},$$
where the sum is over all colorings of~$G$ and the product is over
all vertices $v$ of~$G$.

The reason this observation is significant is that the invariant
$\X_G$ so defined turns out to be
a generalization of the \itc{chromatic polynomial~$\chi_G(n)$} of~$G$.
(This therefore opens up the possibility of applying
the well-established theory of graph colorings to the solution
of the Poset Chain Conjecture.)
More precisely,
recall that~$\chi_G(n)$ is by definition the number of colorings
of~$G$ such that every vertex is assigned a color between 1 and~$n$
inclusive.  From this we see that if we set $n$ of the indeterminates~$x_i$
equal to~1 and the rest equal to~0, $\X_G$ becomes~$\chi_G(n)$.
This specialization procedure (of setting $n$ of the indeterminates
equal to~1 and the rest equal to~0) is an important one,
so we introduce the notation $g(1^n)$ to represent the result of
so specializing the symmetric function~$g$.
We will come back to this specialization later.

Let me now break off this account of this line of research temporarily
in order to describe the other line of research that was alluded to earlier.
This second line of research is the mathematical theory of \itc{juggling.}
Apart from the fact that many mathematicians are amateur jugglers,
the recreational flavor of juggling gives its mathematical study
great intuitive appeal.  In the delightful paper~[BEGW], the authors
show how some intriguing and nontrivial mathematics arises from the
analysis of various juggling patterns.  In particular, they are led
to consider a poset invariant called the \itc{binomial drop
polynomial}~[B-G] (which turns out to be
equivalent to the \itc{factorial polynomial} of~[GJW1]).

In the course of studying the binomial drop polynomial,
Chung and Graham~[CG1] had the idea of generalizing it
to an invariant of a directed graph.
More precisely, they made the simple observation that a poset~$P$
may be thought of as a directed graph~$D(P)$ whose vertices
are the elements of~$P$ and in which $u\to v$ is an edge
if and only if $u<v$ in~$P$.
They were then led to generalize the binomial drop polynomial
of a poset to an invariant of an arbitrary directed graph that
they call the \itc{cover polynomial.}

Since the cover polynomial is central to this thesis,
we digress momentarily to give its definition.
A \itc{path-cycle cover} of a directed graph~$D$ is a spanning
subgraph of~$D$ each of whose connected components is either a
directed path or a directed cycle.  (An isolated vertex is
considered to be a path of zero length and an isolated loop
is considered to be a cycle of length one.)  The cover polynomial
$C(D;i,j)$ is then defined to be
$$C(D;i,j) \defeq \sum_S i^{\underline{p(S)}} j^{c(S)},$$
where the sum is over all path-cycle covers~$S$ of~$D$,
$p(S)$ is the number of connected components of~$S$ that are paths,
$c(S)$ is the number of connected components of~$S$ that are cycles,
and the underline indicates the \itc{lower factorial,} i.e.,
$$i^{\underline k} \defeq i(i-1)\cdots (i-k+1).$$
(This rather mysterious definition is motivated by deletion-contraction
considerations; the interested reader is referred to~[CG2] for more
explanation.)

Returning now to our main account, we remark that based on
what we have said so far, the reader would probably not expect any
connection between the mathematics of immanants and the mathematics
of juggling, other than the fact that posets and graphs enter
the picture in both cases.  The surprise, however, is that for
any poset~$P$, the polynomial (in~$n$) $\X_P(1^n)$ is the
factorial polynomial of~$P$!

This unexpected connection suggests that we ought to stop for
a moment and
clarify the precise relationships between all the algebraic invariants
mentioned so far.  First, in the realm of posets, we have the symmetric
function invariant~$\X_P$.  This becomes, via the specialization
$\X_P \mapsto \X_P(1^n)$, the factorial polynomial or binomial
drop polynomial of~$P$.
The invariant $\X_P$ generalizes to a graph invariant~$\X_G$,
and $\X_G(1^n)$ is the chromatic polynomial of~$G$.
(It follows that the chromatic polynomial is a generalization of
the factorial polynomial; this has actually been known for a long
time~[GJW2].)
We also have a generalization of the factorial polynomial
to a directed graph invariant: the cover polynomial.
It is therefore irresistible to ask (and Chung and Graham
in fact do ask this) if there is a symmetric
function invariant of a directed graph, analogous to~$\X_G$,
which turns into the cover polynomial upon applying the specialization
$g\mapsto g(1^n)$.
This question brings us (finally!) to the main topic of this thesis.

In Chapter~1, among other things, a symmetric function
generalization (the \itc{path-cycle symmetric function})
of the cover polynomial is defined,
following a suggestion of Stanley.
More precisely, the path-cycle symmetric function of a directed graph~$D$
is defined by
$$\Xi_D(x;y) \defeq \sum_S \tilde m_{\pi(S)}(x) p_{\sigma(S)}(y),$$
where the sum is as before over all path-cycle covers~$S$ of~$D$,
$x$ and~$y$ are (countably infinite) sets of commuting independent
indeterminates, $\pi(S)$ denotes the integer partition consisting of
the lengths of the directed paths in~$S$, $\sigma(S)$ denotes the
integer partition consisting of the lengths of the directed cycles
of~$S$, $\tilde m$ denotes the augmented monomial symmetric functions
(the same as the usual monomial symmetric functions except with a 
constant factor), and $p$ denotes the power sum symmetric functions.
(We shall presently give a few more words of motivation for this
unwieldy-looking definition.)

It is not difficult to show that
the path-cycle symmetric function has the desired property of
becoming the cover polynomial under the specialization $g\mapsto g(1^n)$.
In Chapter~2, which forms the bulk of this thesis,
the properties of the path-cycle symmetric function are investigated,
and a number of striking results are obtained.
%
%
In section~1, we obtain analogues of some of the basic theorems
about~$\X_G$.
In section~3, we obtain generalizations of some of the basic results of
rook theory, such as the M\"obius inversion formula of~[GJW4], and the
fundamental inclusion-exclusion formula for rooks.  Here already we obtain
our first surprise.  In the course of generalizing the inclusion-exclusion
formula for rooks, we are led to define a certain (vector space)
basis for symmetric functions that does not seem to have been studied
before but which seems to be an important object.
It has some surprising connections with monomial
symmetric functions and fundamental quasi-symmetric functions
(for example, any symmetric function that is a nonnegative combination
of fundamental quasi-symmetric functions also has nonnegative coefficients
when expanded in terms of this new symmetric function basis),
and it also generalizes the polynomial basis
$$\biggl\{{x+k \choose n}\biggr\}_{k=0,1,\ldots,n}$$
that arises in many contexts in combinatorics.

Perhaps the most surprising result, though, appears in section~2,
which gives a ``reciprocity'' formula relating the
path-cycle symmetric functions of complementary digraphs.
(Two digraphs are complementary if the edges of one are precisely
the non-edges of the other.)  The way this formula was discovered is
instructive: I was trying to find an analogue of [St2, Corollary~2.7],
which is a theorem about~$\omega\X_G$.  Here $\omega$ is an involution
on the space of symmetric functions that arises naturally in many
contexts.  It is therefore natural to consider what $\omega$ does to
the path-cycle symmetric function.  I did not actually succeed in
proving the result I was hoping for, but ended up proving
the reciprocity formula instead.  I was then amazed to discover
afterwards that the reciprocity formula has several beautiful
corollaries.  For example, when specialized to the context of
the cover polynomial, it provides a remarkably simple formula
relating the cover polynomials of complementary digraphs:
$$C(D;i,j) = (-1)^d C(D';-i-j,j),$$
where $d$ is the number of vertices of~$D$.
(Chung and Graham had conjectured that some such relation might
exist, but did not even have a conjectural formula.  In fact,
the reciprocity formula for cover polynomials is not an easy
one to guess, and it seems that generalizing to the symmetric
function context actually makes it easier to find the
pattern---although Ira Gessel did independently arrive
at the reciprocity formula for cover polynomials without
resorting to symmetric functions.)
Further specialization of the reciprocity formula to the context
of factorial polynomials gives a formula relating factorial
polynomials of complementary boards that is much simpler and
more elegant than previously known formulas.  Also, a different
specialization of the reciprocity formula gives a result
that Stanley and Stembridge have proved in the course of their study
of immanants.


The remainder of the thesis consists of a miscellaneous collection
of results about the path-cycle symmetric function
(scattered throughout Chapter~2) and about~$\X_G$ (in Chapter~3).
For example, $\X_G$ is shown to be reconstructible from the list
of vertex-deleted subgraphs of~$G$, and Chung and Graham's concept
of $G$-descents is shown to have an application to the study
of~$\X_G$.  Instead of giving further details about these results,
however, I wish to conclude by addressing the
question: where do we go from here?

We have seen that $\X_G$ and the path-cycle symmetric function
are natural objects to study (if one is familiar with the relevant
background), and the attractive
results obtained in the course of investigating
them also provide some a posteriori justification for their study.
But once we have derived the main properties of the path-cycle
symmetric function, what is there left to do?  Is this a dead end?

The title of this thesis suggests a possible direction for further
research.  We know now that the chromatic polynomial and the cover
polynomial have interesting symmetric function generalizations.
What about other graph polynomials, or more generally, other
combinatorial polynomials?  Do they have natural symmetric function
generalizations?

To answer this question, let us re-examine how the chromatic
polynomial and the cover polynomial are generalized.
In both cases, it turns out that the method of generalization
is based on interpreting the polynomial as counting the \itc{total}
number of colorings of a certain kind.  The generalization is then
obtained by enumerating the same set of colorings, but keeping track
of the number of times \itc{each color} is used.  (This is the
motivation for the unusual-looking definition of~$\Xi_D$ given above.)

The problem with this idea is that most combinatorial polynomials
do not have any obvious interpretation in terms of colorings.
However, it is possible to modify the above idea slightly so that
it applies to a wider class of algebraic invariants.  I will not
go into details here, since I hope to do so elsewhere, but the
essential idea is to generalize a polynomial by first interpreting
it as enumerating the total number of \itc{partitions of a set}
of a certain kind, and then obtaining a generalization by
enumerating the same set of partitions but keeping track of the
\itc{type} of the partition.  There are already some indications
that this approach ought to be fruitful; for example, Doubilet's
theory of symmetric functions [Dou] is based on this idea (and
this probably accounts for why his interpretations of the
change-of-basis coefficients between the various standard
symmetric function bases have proved to be more useful in our
work than the formulas in [Mac, Table~1, p.~56] have).
The theory of polynomials of \itc{binomial type}~[R-R]
also generalizes nicely using this idea; the binomial type
identity generalizes to the coproduct in the Hopf algebra
of symmetric functions.  In short, there are many promising
directions for further study.

\vfill\eject

{\nopagenumbers
\hbox{ }
\vfill\eject}

\def\PsfigVersion{1.10}
\def\setDriver{\DvipsDriver} 
\ifx\undefined\psfig\else\endinput\fi
%

\let\LaTeXAtSign=\@
\let\@=\relax
\edef\psfigRestoreAt{\catcode`\@=\number\catcode`@\relax}
\catcode`\@=11\relax
\newwrite\@unused
\def\ps@typeout#1{{\let\protect\string\immediate\write\@unused{#1}}}

\def\DvipsDriver{
	\ps@typeout{psfig/tex \PsfigVersion -dvips}
\def\PsfigSpecials{\DvipsSpecials} 	\def\ps@dir{/}
\def\ps@predir{} }
\def\OzTeXDriver{
	\ps@typeout{psfig/tex \PsfigVersion -oztex}
	\def\PsfigSpecials{\OzTeXSpecials}
	\def\ps@dir{:}
	\def\ps@predir{:}
	\catcode`\^^J=5
}


\def\figurepath{./:}

\def\DoPaths#1{\expandafter\EachPath#1\stoplist}
\def\leer{}
\def\EachPath#1:#2\stoplist{
  \ExistsFile{#1}{\SearchedFile}
  \ifx#2\leer
  \else
    \expandafter\EachPath#2\stoplist
  \fi}
%
%
\def\ps@dir{/}
\def\ExistsFile#1#2{%
   \openin1=\ps@predir#1\ps@dir#2
   \ifeof1
       \closein1
   \else
       \closein1
        \ifx\ps@founddir\leer
           \edef\ps@founddir{#1}
        \fi
   \fi}
%
%
\def\get@dir#1{%
  \def\ps@founddir{}
  \def\SearchedFile{#1}
  \DoPaths\figurepath
}

%
%
\def\@nnil{\@nil}
\def\@empty{}
\def\@psdonoop#1\@@#2#3{}
\def\@psdo#1:=#2\do#3{\edef\@psdotmp{#2}\ifx\@psdotmp\@empty \else
    \expandafter\@psdoloop#2,\@nil,\@nil\@@#1{#3}\fi}
\def\@psdoloop#1,#2,#3\@@#4#5{\def#4{#1}\ifx #4\@nnil \else
       #5\def#4{#2}\ifx #4\@nnil \else#5\@ipsdoloop #3\@@#4{#5}\fi\fi}
\def\@ipsdoloop#1,#2\@@#3#4{\def#3{#1}\ifx #3\@nnil 
       \let\@nextwhile=\@psdonoop \else
      #4\relax\let\@nextwhile=\@ipsdoloop\fi\@nextwhile#2\@@#3{#4}}
\def\@tpsdo#1:=#2\do#3{\xdef\@psdotmp{#2}\ifx\@psdotmp\@empty \else
    \@tpsdoloop#2\@nil\@nil\@@#1{#3}\fi}
\def\@tpsdoloop#1#2\@@#3#4{\def#3{#1}\ifx #3\@nnil 
       \let\@nextwhile=\@psdonoop \else
      #4\relax\let\@nextwhile=\@tpsdoloop\fi\@nextwhile#2\@@#3{#4}}
%
\ifx\undefined\fbox
\newdimen\fboxrule
\newdimen\fboxsep
\newdimen\ps@tempdima
\newbox\ps@tempboxa
\fboxsep = 3pt
\fboxrule = .4pt
\long\def\fbox#1{\leavevmode\setbox\ps@tempboxa\hbox{#1}\ps@tempdima\fboxrule
    \advance\ps@tempdima \fboxsep \advance\ps@tempdima \dp\ps@tempboxa
   \hbox{\lower \ps@tempdima\hbox
  {\vbox{\hrule height \fboxrule
          \hbox{\vrule width \fboxrule \hskip\fboxsep
          \vbox{\vskip\fboxsep \box\ps@tempboxa\vskip\fboxsep}\hskip 
                 \fboxsep\vrule width \fboxrule}
                 \hrule height \fboxrule}}}}
\fi
%
%
\newread\ps@stream
\newif\ifnot@eof       
\newif\if@noisy        
\newif\if@atend        
\newif\if@psfile       
%
%
{\catcode`\%=12\global\gdef\epsf@start{
\def\epsf@PS{PS}
\def\epsf@getbb#1{%
%
%
\openin\ps@stream=\ps@predir#1
\ifeof\ps@stream\ps@typeout{Error, File #1 not found}\else
%
%
   {\not@eoftrue \chardef\other=12
    \def\do##1{\catcode`##1=\other}\dospecials \catcode`\ =10
    \loop
       \if@psfile
	  \read\ps@stream to \epsf@fileline
       \else{
	  \obeyspaces
          \read\ps@stream to \epsf@tmp\global\let\epsf@fileline\epsf@tmp}
       \fi
       \ifeof\ps@stream\not@eoffalse\else
%
%
       \if@psfile\else
       \expandafter\epsf@test\epsf@fileline:. \\%
       \fi
%
%
          \expandafter\epsf@aux\epsf@fileline:. \\%
       \fi
   \ifnot@eof\repeat
   }\closein\ps@stream\fi}%
%
%
\long\def\epsf@test#1#2#3:#4\\{\def\epsf@testit{#1#2}
			\ifx\epsf@testit\epsf@start\else
\ps@typeout{Warning! File does not start with `\epsf@start'.  It may not be a PostScript file.}
			\fi
			\@psfiletrue} 
%
%
{\catcode`\%=12\global\let\epsf@percent=
%
%
%
\long\def\epsf@aux#1#2:#3\\{\ifx#1\epsf@percent
   \def\epsf@testit{#2}\ifx\epsf@testit\epsf@bblit
	\@atendfalse
        \epsf@atend #3 . \\%
	\if@atend	
	   \if@verbose{
		\ps@typeout{psfig: found `(atend)'; continuing search}
	   }\fi
        \else
        \epsf@grab #3 . . . \\%
        \not@eoffalse
        \global\no@bbfalse
        \fi
   \fi\fi}%
%
%
\def\epsf@grab #1 #2 #3 #4 #5\\{%
   \global\def\epsf@llx{#1}\ifx\epsf@llx\empty
      \epsf@grab #2 #3 #4 #5 .\\\else
   \global\def\epsf@lly{#2}%
   \global\def\epsf@urx{#3}\global\def\epsf@ury{#4}\fi}%
%
%
\def\epsf@atendlit{(atend)} 
\def\epsf@atend #1 #2 #3\\{%
   \def\epsf@tmp{#1}\ifx\epsf@tmp\empty
      \epsf@atend #2 #3 .\\\else
   \ifx\epsf@tmp\epsf@atendlit\@atendtrue\fi\fi}


\chardef\psletter = 11 
\chardef\other = 12

\newif \ifdebug 
\newif\ifc@mpute 
\c@mputetrue 

\let\then = \relax
\def\r@dian{pt }
\let\r@dians = \r@dian
\let\dimensionless@nit = \r@dian
\let\dimensionless@nits = \dimensionless@nit
\def\internal@nit{sp }
\let\internal@nits = \internal@nit
\newif\ifstillc@nverging
\def \Mess@ge #1{\ifdebug \then \message {#1} \fi}

{ 
	\catcode `\@ = \psletter
	\gdef \nodimen {\expandafter \n@dimen \the \dimen}
	\gdef \term #1 #2 #3%
	       {\edef \t@ {\the #1}
		\edef \t@@ {\expandafter \n@dimen \the #2\r@dian}%
		\t@rm {\t@} {\t@@} {#3}%
	       }
	\gdef \t@rm #1 #2 #3%
	       {{%
		\count 0 = 0
		\dimen 0 = 1 \dimensionless@nit
		\dimen 2 = #2\relax
		\Mess@ge {Calculating term #1 of \nodimen 2}%
		\loop
		\ifnum	\count 0 < #1
		\then	\advance \count 0 by 1
			\Mess@ge {Iteration \the \count 0 \space}%
			\Multiply \dimen 0 by {\dimen 2}%
			\Mess@ge {After multiplication, term = \nodimen 0}%
			\Divide \dimen 0 by {\count 0}%
			\Mess@ge {After division, term = \nodimen 0}%
		\repeat
		\Mess@ge {Final value for term #1 of 
				\nodimen 2 \space is \nodimen 0}%
		\xdef \Term {#3 = \nodimen 0 \r@dians}%
		\aftergroup \Term
	       }}
	\catcode `\p = \other
	\catcode `\t = \other
	\gdef \n@dimen #1pt{#1} 
}

\def \Divide #1by #2{\divide #1 by #2} 

\def \Multiply #1by #2
       {{
	\count 0 = #1\relax
	\count 2 = #2\relax
	\count 4 = 65536
	\Mess@ge {Before scaling, count 0 = \the \count 0 \space and
			count 2 = \the \count 2}%
	\ifnum	\count 0 > 32767 
	\then	\divide \count 0 by 4
		\divide \count 4 by 4
	\else	\ifnum	\count 0 < -32767
		\then	\divide \count 0 by 4
			\divide \count 4 by 4
		\else
		\fi
	\fi
	\ifnum	\count 2 > 32767 
	\then	\divide \count 2 by 4
		\divide \count 4 by 4
	\else	\ifnum	\count 2 < -32767
		\then	\divide \count 2 by 4
			\divide \count 4 by 4
		\else
		\fi
	\fi
	\multiply \count 0 by \count 2
	\divide \count 0 by \count 4
	\xdef \product {#1 = \the \count 0 \internal@nits}%
	\aftergroup \product
       }}

\def\r@duce{\ifdim\dimen0 > 90\r@dian \then   
		\multiply\dimen0 by -1
		\advance\dimen0 by 180\r@dian
		\r@duce
	    \else \ifdim\dimen0 < -90\r@dian \then  
		\advance\dimen0 by 360\r@dian
		\r@duce
		\fi
	    \fi}

\def\Sine#1%
       {{%
	\dimen 0 = #1 \r@dian
	\r@duce
	\ifdim\dimen0 = -90\r@dian \then
	   \dimen4 = -1\r@dian
	   \c@mputefalse
	\fi
	\ifdim\dimen0 = 90\r@dian \then
	   \dimen4 = 1\r@dian
	   \c@mputefalse
	\fi
	\ifdim\dimen0 = 0\r@dian \then
	   \dimen4 = 0\r@dian
	   \c@mputefalse
	\fi
	\ifc@mpute \then
		\divide\dimen0 by 180
		\dimen0=3.141592654\dimen0
		\dimen 2 = 3.1415926535897963\r@dian 
		\divide\dimen 2 by 2 
		\Mess@ge {Sin: calculating Sin of \nodimen 0}%
		\count 0 = 1 
		\dimen 2 = 1 \r@dian 
		\dimen 4 = 0 \r@dian 
		\loop
			\ifnum	\dimen 2 = 0 
			\then	\stillc@nvergingfalse 
			\else	\stillc@nvergingtrue
			\fi
			\ifstillc@nverging 
			\then	\term {\count 0} {\dimen 0} {\dimen 2}%
				\advance \count 0 by 2
				\count 2 = \count 0
				\divide \count 2 by 2
				\ifodd	\count 2 
				\then	\advance \dimen 4 by \dimen 2
				\else	\advance \dimen 4 by -\dimen 2
				\fi
		\repeat
	\fi		
			\xdef \sine {\nodimen 4}%
       }}

\def\Cosine#1{\ifx\sine\UnDefined\edef\Savesine{\relax}\else
		             \edef\Savesine{\sine}\fi
	{\dimen0=#1\r@dian\advance\dimen0 by 90\r@dian
	 \Sine{\nodimen 0}
	 \xdef\cosine{\sine}
	 \xdef\sine{\Savesine}}}	      

\def\psdraft{
	\def\@psdraft{0}
}
\def\psfull{
	\def\@psdraft{100}
}

\psfull

\newif\if@scalefirst
\def\psscalefirst{\@scalefirsttrue}
\def\psrotatefirst{\@scalefirstfalse}
\psrotatefirst

\newif\if@draftbox
\def\psnodraftbox{
	\@draftboxfalse
}
\def\psdraftbox{
	\@draftboxtrue
}
\@draftboxtrue

\newif\if@prologfile
\newif\if@postlogfile
\def\pssilent{
	\@noisyfalse
}
\def\psnoisy{
	\@noisytrue
}
\psnoisy
\newif\if@bbllx
\newif\if@bblly
\newif\if@bburx
\newif\if@bbury
\newif\if@height
\newif\if@width
\newif\if@rheight
\newif\if@rwidth
\newif\if@angle
\newif\if@clip
\newif\if@verbose
\def\@p@@sclip#1{\@cliptrue}
\newif\if@decmpr
\def\@p@@sfigure#1{\def\@p@sfile{null}\def\@p@sbbfile{null}\@decmprfalse
   \openin1=\ps@predir#1
   \ifeof1
	\closein1
	\get@dir{#1}
	\ifx\ps@founddir\leer
		\openin1=\ps@predir#1.bb
		\ifeof1
			\closein1
			\get@dir{#1.bb}
			\ifx\ps@founddir\leer
				\ps@typeout{Can't find #1 in \figurepath}
			\else
				\@decmprtrue
				\def\@p@sfile{\ps@founddir\ps@dir#1}
				\def\@p@sbbfile{\ps@founddir\ps@dir#1.bb}
			\fi
		\else
			\closein1
			\@decmprtrue
			\def\@p@sfile{#1}
			\def\@p@sbbfile{#1.bb}
		\fi
	\else
		\def\@p@sfile{\ps@founddir\ps@dir#1}
		\def\@p@sbbfile{\ps@founddir\ps@dir#1}
	\fi
   \else
	\closein1
	\def\@p@sfile{#1}
	\def\@p@sbbfile{#1}
   \fi
}
\def\@p@@sfile#1{\@p@@sfigure{#1}}
\def\@p@@sbbllx#1{
		\@bbllxtrue
		\dimen100=#1
		\edef\@p@sbbllx{\number\dimen100}
}
\def\@p@@sbblly#1{
		\@bbllytrue
		\dimen100=#1
		\edef\@p@sbblly{\number\dimen100}
}
\def\@p@@sbburx#1{
		\@bburxtrue
		\dimen100=#1
		\edef\@p@sbburx{\number\dimen100}
}
\def\@p@@sbbury#1{
		\@bburytrue
		\dimen100=#1
		\edef\@p@sbbury{\number\dimen100}
}
\def\@p@@sheight#1{
		\@heighttrue
		\dimen100=#1
   		\edef\@p@sheight{\number\dimen100}
}
\def\@p@@swidth#1{
		\@widthtrue
		\dimen100=#1
		\edef\@p@swidth{\number\dimen100}
}
\def\@p@@srheight#1{
		\@rheighttrue
		\dimen100=#1
		\edef\@p@srheight{\number\dimen100}
}
\def\@p@@srwidth#1{
		\@rwidthtrue
		\dimen100=#1
		\edef\@p@srwidth{\number\dimen100}
}
\def\@p@@sangle#1{
		\@angletrue
		\edef\@p@sangle{#1} 
}
\def\@p@@ssilent#1{ 
		\@verbosefalse
}
\def\@p@@sprolog#1{\@prologfiletrue\def\@prologfileval{#1}}
\def\@p@@spostlog#1{\@postlogfiletrue\def\@postlogfileval{#1}}
\def\@cs@name#1{\csname #1\endcsname}
\def\@setparms#1=#2,{\@cs@name{@p@@s#1}{#2}}
%
%
\def\ps@init@parms{
		\@bbllxfalse \@bbllyfalse
		\@bburxfalse \@bburyfalse
		\@heightfalse \@widthfalse
		\@rheightfalse \@rwidthfalse
		\def\@p@sbbllx{}\def\@p@sbblly{}
		\def\@p@sbburx{}\def\@p@sbbury{}
		\def\@p@sheight{}\def\@p@swidth{}
		\def\@p@srheight{}\def\@p@srwidth{}
		\def\@p@sangle{0}
		\def\@p@sfile{} \def\@p@sbbfile{}
		\def\@p@scost{10}
		\def\@sc{}
		\@prologfilefalse
		\@postlogfilefalse
		\@clipfalse
		\if@noisy
			\@verbosetrue
		\else
			\@verbosefalse
		\fi
}
%
%
\def\parse@ps@parms#1{
	 	\@psdo\@psfiga:=#1\do
		   {\expandafter\@setparms\@psfiga,}}
%
%
\newif\ifno@bb
\def\bb@missing{
	\if@verbose{
		\ps@typeout{psfig: searching \@p@sbbfile \space  for bounding box}
	}\fi
	\no@bbtrue
	\epsf@getbb{\@p@sbbfile}
        \ifno@bb \else \bb@cull\epsf@llx\epsf@lly\epsf@urx\epsf@ury\fi
}	
\def\bb@cull#1#2#3#4{
	\dimen100=#1 bp\edef\@p@sbbllx{\number\dimen100}
	\dimen100=#2 bp\edef\@p@sbblly{\number\dimen100}
	\dimen100=#3 bp\edef\@p@sbburx{\number\dimen100}
	\dimen100=#4 bp\edef\@p@sbbury{\number\dimen100}
	\no@bbfalse
}
\newdimen\p@intvaluex
\newdimen\p@intvaluey
\def\rotate@#1#2{{\dimen0=#1 sp\dimen1=#2 sp
		  \global\p@intvaluex=\cosine\dimen0
		  \dimen3=\sine\dimen1
		  \global\advance\p@intvaluex by -\dimen3
		  \global\p@intvaluey=\sine\dimen0
		  \dimen3=\cosine\dimen1
		  \global\advance\p@intvaluey by \dimen3
		  }}
\def\compute@bb{
		\no@bbfalse
		\if@bbllx \else \no@bbtrue \fi
		\if@bblly \else \no@bbtrue \fi
		\if@bburx \else \no@bbtrue \fi
		\if@bbury \else \no@bbtrue \fi
		\ifno@bb \bb@missing \fi
		\ifno@bb \ps@typeout{FATAL ERROR: no bb supplied or found}
			\no-bb-error
		\fi
		%
%
		\count203=\@p@sbburx
		\count204=\@p@sbbury
		\advance\count203 by -\@p@sbbllx
		\advance\count204 by -\@p@sbblly
		\edef\ps@bbw{\number\count203}
		\edef\ps@bbh{\number\count204}
		\if@angle 
			\Sine{\@p@sangle}\Cosine{\@p@sangle}
	        	{\dimen100=\maxdimen\xdef\r@p@sbbllx{\number\dimen100}
					    \xdef\r@p@sbblly{\number\dimen100}
			                    \xdef\r@p@sbburx{-\number\dimen100}
					    \xdef\r@p@sbbury{-\number\dimen100}}
%
                        \def\minmaxtest{
			   \ifnum\number\p@intvaluex<\r@p@sbbllx
			      \xdef\r@p@sbbllx{\number\p@intvaluex}\fi
			   \ifnum\number\p@intvaluex>\r@p@sbburx
			      \xdef\r@p@sbburx{\number\p@intvaluex}\fi
			   \ifnum\number\p@intvaluey<\r@p@sbblly
			      \xdef\r@p@sbblly{\number\p@intvaluey}\fi
			   \ifnum\number\p@intvaluey>\r@p@sbbury
			      \xdef\r@p@sbbury{\number\p@intvaluey}\fi
			   }
			\rotate@{\@p@sbbllx}{\@p@sbblly}
			\minmaxtest
			\rotate@{\@p@sbbllx}{\@p@sbbury}
			\minmaxtest
			\rotate@{\@p@sbburx}{\@p@sbblly}
			\minmaxtest
			\rotate@{\@p@sbburx}{\@p@sbbury}
			\minmaxtest
			\edef\@p@sbbllx{\r@p@sbbllx}\edef\@p@sbblly{\r@p@sbblly}
			\edef\@p@sbburx{\r@p@sbburx}\edef\@p@sbbury{\r@p@sbbury}
		\fi
		\count203=\@p@sbburx
		\count204=\@p@sbbury
		\advance\count203 by -\@p@sbbllx
		\advance\count204 by -\@p@sbblly
		\edef\@bbw{\number\count203}
		\edef\@bbh{\number\count204}
}
%
%
\def\in@hundreds#1#2#3{\count240=#2 \count241=#3
		     \count100=\count240	
		     \divide\count100 by \count241
		     \count101=\count100
		     \multiply\count101 by \count241
		     \advance\count240 by -\count101
		     \multiply\count240 by 10
		     \count101=\count240	
		     \divide\count101 by \count241
		     \count102=\count101
		     \multiply\count102 by \count241
		     \advance\count240 by -\count102
		     \multiply\count240 by 10
		     \count102=\count240	
		     \divide\count102 by \count241
		     \count200=#1\count205=0
		     \count201=\count200
			\multiply\count201 by \count100
		 	\advance\count205 by \count201
		     \count201=\count200
			\divide\count201 by 10
			\multiply\count201 by \count101
			\advance\count205 by \count201
		     \count201=\count200
			\divide\count201 by 100
			\multiply\count201 by \count102
			\advance\count205 by \count201
		     \edef\@result{\number\count205}
}
\def\compute@wfromh{
		\in@hundreds{\@p@sheight}{\@bbw}{\@bbh}
		\edef\@p@swidth{\@result}
}
\def\compute@hfromw{
	        \in@hundreds{\@p@swidth}{\@bbh}{\@bbw}
		\edef\@p@sheight{\@result}
}
\def\compute@handw{
		\if@height 
			\if@width
			\else
				\compute@wfromh
			\fi
		\else 
			\if@width
				\compute@hfromw
			\else
				\edef\@p@sheight{\@bbh}
				\edef\@p@swidth{\@bbw}
			\fi
		\fi
}
\def\compute@resv{
		\if@rheight \else \edef\@p@srheight{\@p@sheight} \fi
		\if@rwidth \else \edef\@p@srwidth{\@p@swidth} \fi
}
%
\def\compute@sizes{
	\compute@bb
	\if@scalefirst\if@angle
	\if@width
	   \in@hundreds{\@p@swidth}{\@bbw}{\ps@bbw}
	   \edef\@p@swidth{\@result}
	\fi
	\if@height
	   \in@hundreds{\@p@sheight}{\@bbh}{\ps@bbh}
	   \edef\@p@sheight{\@result}
	\fi
	\fi\fi
	\compute@handw
	\compute@resv}
\def\OzTeXSpecials{
	\special{empty.ps /@isp {true} def}
	\special{empty.ps \@p@swidth \space \@p@sheight \space
			\@p@sbbllx \space \@p@sbblly \space
			\@p@sbburx \space \@p@sbbury \space
			startTexFig \space }
	\if@clip{
		\if@verbose{
			\ps@typeout{(clip)}
		}\fi
		\special{empty.ps doclip \space }
	}\fi
	\if@angle{
		\if@verbose{
			\ps@typeout{(rotate)}
		}\fi
		\special {empty.ps \@p@sangle \space rotate \space} 
	}\fi
	\if@prologfile
	    \special{\@prologfileval \space } \fi
	\if@decmpr{
		\if@verbose{
			\ps@typeout{psfig: Compression not available
			in OzTeX version \space }
		}\fi
	}\else{
		\if@verbose{
			\ps@typeout{psfig: including \@p@sfile \space }
		}\fi
		\special{epsf=\ps@predir\@p@sfile \space }
	}\fi
	\if@postlogfile
	    \special{\@postlogfileval \space } \fi
	\special{empty.ps /@isp {false} def}
}
\def\DvipsSpecials{
	\special{ps::[begin] 	\@p@swidth \space \@p@sheight \space
			\@p@sbbllx \space \@p@sbblly \space
			\@p@sbburx \space \@p@sbbury \space
			startTexFig \space }
	\if@clip{
		\if@verbose{
			\ps@typeout{(clip)}
		}\fi
		\special{ps:: doclip \space }
	}\fi
	\if@angle
		\if@verbose{
			\ps@typeout{(clip)}
		}\fi
		\special {ps:: \@p@sangle \space rotate \space} 
	\fi
	\if@prologfile
	    \special{ps: plotfile \@prologfileval \space } \fi
	\if@decmpr{
		\if@verbose{
			\ps@typeout{psfig: including \@p@sfile.Z \space }
		}\fi
		\special{ps: plotfile "`zcat \@p@sfile.Z" \space }
	}\else{
		\if@verbose{
			\ps@typeout{psfig: including \@p@sfile \space }
		}\fi
		\special{ps: plotfile \@p@sfile \space }
	}\fi
	\if@postlogfile
	    \special{ps: plotfile \@postlogfileval \space } \fi
	\special{ps::[end] endTexFig \space }
}
%
%
\def\psfig#1{\vbox {
	%
	\ps@init@parms
	\parse@ps@parms{#1}
	\compute@sizes
	\ifnum\@p@scost<\@psdraft{
		\PsfigSpecials 
		\vbox to \@p@srheight sp{
			\hbox to \@p@srwidth sp{
				\hss
			}
		\vss
		}
	}\else{
		\if@draftbox{		
			\hbox{\fbox{\vbox to \@p@srheight sp{
			\vss
			\hbox to \@p@srwidth sp{ \hss 
			 \hss }
			\vss
			}}}
		}\else{
			\vbox to \@p@srheight sp{
			\vss
			\hbox to \@p@srwidth sp{\hss}
			\vss
			}
		}\fi

	}\fi
}}
\psfigRestoreAt
\setDriver
\let\@=\LaTeXAtSign

\hfuzz=3pt

\newcount\chapno \chapno=0
\def\chap#1\par{\global\advance\chapno by 1\vfill\eject
  \vbox to 1.5truein{\vfill}\leftline{\titlefont Chapter \the\chapno}\vbox
        to 1truein{\vfill}\leftline{\titlefont #1}\vbox
        to 1.5truein{\vfill}\noindent \global\secno=0}

\newcount\thmno \thmno=1
\long\def\theorem#1#2{\edef#1{Theorem~\the\thmno}\noindent
    {\sc Theorem \the\thmno.\enspace}{\it #2}\endgraf
    \penalty55\global\advance\thmno by 1}
\newcount\lemmano \lemmano=1
\long\def\lemma#1#2{\edef#1{Lemma~\the\lemmano}\noindent
    {\sc Lemma \the\lemmano.\enspace}{\it #2}\endgraf
    \penalty55\global\advance\lemmano by 1}
\newcount\propno \propno=1
\long\def\proposition#1#2{\edef#1{Proposition~\the\propno}\noindent
    {\sc Proposition \the\propno.\enspace}{\it #2}\endgraf
    \penalty55\global\advance\propno by 1}
\newcount\corolno \corolno=1
\long\def\corollary#1#2{\edef#1{Corollary~\the\corolno}\noindent
    {\sc Corollary \the\corolno.\enspace}{\it #2}\endgraf
    \penalty55\global\advance\corolno by 1}
\newcount\conjno \conjno=1
\long\def\conjecture#1#2{\edef#1{Conjecture~\the\conjno}\noindent
    {\sc Conjecture \the\conjno.\enspace}{\it #2}\endgraf
    \global\advance\conjno by 1}

\newcount\secno \secno=0
\def\sec#1\par{\global\advance\secno by 1\bigskip
    \noindent\centerline{\sc\the\secno. #1}\par\nobreak\noindent}

\def\definition{\noindent{\sc Definition.\enspace}}
\def\proof{\noindent{\it Proof.\enspace}}
\def\qed{~\vrule width4pt height6pt depth1pt}

\font\noterm=cmr9
\font\noteit=cmti9
\catcode`\@=11
\def\footnote#1{\edef\@sf{\spacefactor\the\spacefactor}#1\@sf
      \insert\footins\bgroup\noterm
      \interlinepenalty100 \let\par=\endgraf
        \leftskip=\z@skip \rightskip=\z@skip
        \splittopskip=10pt plus 1pt minus 1pt \floatingpenalty=20000
        \medskip\item{#1}\bgroup\strut\aftergroup\@foot\let\next}
\catcode`@=12
\dimen\footins=30pc 
\newdimen\notespc \notespc=11pt
\newcount\notenum \notenum=0
\newdimen\currentspc
\def\note#1{\currentspc=\the\baselineskip
      \parskip=0pt\baselineskip=\notespc
      \global\advance\notenum by 1\footnote{$^{\the\notenum}$}\bgroup
      #1\egroup { }\baselineskip=\currentspc \parskip=12pt plus4pt minus4pt}

\parskip=12pt plus4pt minus4pt
\font\bigbold=cmbx12
\font\sc=cmcsc10
\def\itc#1{{\it #1\/}}
\def\X{{\rm X}}
\def\sgn{{\rm sgn}\,}


\chap Preliminaries

In this chapter, we provide some basic definitions and facts.
The absence of a reference does not necessarily imply a claim
to originality, since some references for ``standard'' concepts
and facts have been omitted.

\sec Graphs, digraphs, and symmetric functions

We shall assume that reader is familiar with the basic facts about
set partitions, posets, M\"obius functions, permutations, and so on;
a good reference is~[St3].

Throughout, the unadorned term \itc{graph} will mean a
finite simple labelled undirected graph and the term \itc{digraph}
will mean a finite labelled directed graph without multiple edges
but possibly with loops and bidirected edges.
If $G$ is a graph or a digraph we let $V(G)$
and~$E(G)$ denote its vertex set and edge set respectively.
If $d$ is a positive integer, we use the notation
$[d]$ for the set $\{1,2,\ldots,d\}$.
Note that with our conventions, a digraph~$D$ with $d$~vertices
is equivalent to a subset of $[d]\times[d]$, i.e., a \itc{board.}
(Consider the edge set of~$D$.)  We call this subset the \itc{associated
board,} and conversely given a board we call the corresponding digraph
on~$[d]$ the \itc{associated digraph.}
Since the two representations are equivalent, we shall switch
freely between them, sometimes without warning.

If $G$ is a graph, we say that
a partition~$\sigma$ of~$V(G)$ is \itc{connected} if for each block~$B$
of~$\sigma$ the subgraph induced by~$B$ is connected.  The set of all
connected partitions of~$G$, partially ordered by refinement, forms a
lattice~$L_G$ called the \itc{lattice of contractions} or \itc{bond
lattice} of~$G$.

Our notation for symmetric functions and partitions for the most
part follows that of Macdonald~[Mac], to which the reader is referred for
any facts about symmetric functions that we do not explicitly reference.
We shall always deal with symmetric functions in countably many variables.
If $\tau$ is a set partition or an integer partition, we
write~$\ell(\tau)$ for the number of parts of~$\tau$,
and $|\tau|$ for the sum of the sizes of the parts of~$\tau$.
We define $\sgn\tau$ by
$$\sgn\tau \defeq (-1)^{|\tau| - \ell(\tau)}.$$
We also define
$$r_\tau! \defeq r_1!\,r_2!\cdots,$$
where $r_i$ is the number of parts of~$\tau$ of size~$i$.
We shall denote symmetric functions by a single letter such as~$g$
or by~$g(x)$ if we wish to emphasize that the symmetric function is
in the variables $x=(x_1,x_2,\ldots)$.
Similarly, the ring of symmetric functions will be denoted by~$\Lambda$
or~$\Lambda(x)$.
In addition to the usual symmetric functions $m_\lambda$, $p_\lambda$,
$e_\lambda$, $h_\lambda$, and $s_\lambda$, we shall need the \itc{augmented
monomial symmetric functions~$\tilde m_\lambda$}~[D-K], which are defined by
$$\tilde m_\lambda \defeq r_\lambda!\,m_\lambda.$$
We shall also need the
\itc{forgotten symmetric functions~$f_\lambda$,} which are defined by
$$f_{\lambda} \defeq (\sgn\lambda) \, \omega(\tilde m_\lambda).$$
(Warning: this is one place where we deviate from Macdonald's conventions
and follow Doubilet~[Dou] instead, since [Dou] contains all the results
about the forgotten symmetric functions that we shall need.)
The symbol $\omega$ denotes the usual
involution on symmetric functions that sends $e_\lambda$ to~$h_\lambda$.
If $g(x)$ is a symmetric function, we shall write
$g(-x)$ for the function obtained by negating each variable,
and we shall write $g(1^n)$ for the polynomial in the variable~$n$
obtained by setting $n$ variables equal to one and the rest equal to zero.
We will sometimes use
\itc{set} partitions instead of \itc{integer} partitions in subscripts;
for example, if $\pi$ is a set partition then the expression $p_\pi$ is
to be understood as an abbreviation for~$p_{\rm type(\pi)}$.
We also say that a symmetric function~$g$ is
\itc{$u$-positive} if $\{u_\lambda\}$ is a symmetric function basis and
the expansion of~$g$ in terms of this basis has nonnegative coefficients.

We shall be dealing frequently with functions in two sets of variables,
i.e., elements of $\Lambda(x)\otimes\Lambda(y)$,
so we fix some notation here.  Let $\{x_1,x_2,x_3,\ldots\}$ and
$\{y_1,y_2,y_3,\ldots\}$ be two sets of independent indeterminates.
(Everything commutes with everything else.)  An expression like $g(x;y)$
indicates that $g$ is invariant under any permutation of the $x$~variables
and any permutation of the $y$~variables.  If $g$ is also invariant under
permutations that mix $x$ and~$y$ variables, then we will sometimes write
$g(x,y)$ instead of~$g(x;y)$.  For example, $p_\lambda(x,y)$
indicates the power sum symmetric function in the union of the $x$
and~$y$ variables.  Expressions like $g(x;0)$, $g(x;-y)$ and $g(1^i;1^j)$
have their natural meanings.
The notation $\omega_x g$ will indicate that
for the purposes of applying~$\omega$, $g$ is to be interpreted as a
symmetric function in the $x$ variables with coefficients in the~$y$'s.

\sec $\chi_G$, $\X_G$, $C(D)$, and $\Xi_D$

We now turn to the definitions of the polynomial and symmetric function
invariants that we shall be studying.

Let $G$ be a graph.  A \itc{proper coloring} or simply a \itc{coloring}
of~$G$ is an assignment of positive integers to the vertices of~$G$ such
that adjacent vertices are never assigned the same integer.
The \itc{chromatic polynomial~$\chi_G(i)$} of a graph~$G$ is defined
to be the number of colorings of~$G$ using at most $i$ colors (i.e.,
such that the set of assigned integers is a subset of~$[i]$).
To see that $\chi_G(i)$ is indeed a polynomial in~$i$, let us define
a \itc{stable partition} of~$G$ to be a partition of~$V(G)$
such that no two vertices in the same block are connected by an edge.
Now observe that if we take any stable partition of~$G$ and assign
integers in~$[i]$ to the vertices of~$G$ in such a way that
two vertices are assigned the same integer
if and only if they are in the same block, then we
obtain a coloring of~$G$, and that moreover each coloring of~$G$
can be constructed in one and only one way by this procedure.
Hence
$$\chi_G(i) = \sum_\pi i^{\underline{\ell(\pi)}},$$
where the sum is over all stable partitions~$\pi$ of~$G$.
(We are using the notation $i^{\underline k} = i(i-1)\cdots(i-k+1)$
and $i^{\overline k} = i(i+1)\cdots(i+k-1)$.)
In particular, $\chi_G(i)$ is a polynomial in~$i$, as advertised.

Now let $x_1, x_2, \ldots$ be commuting independent indeterminates.
Stanley~[St2] defines the \itc{chromatic symmetric function~$\X_G$} by
$$\X_G = \X_G(x) \defeq \sum_\pi \tilde m_\pi(x),$$
where the sum is over all stable partitions of~$G$.  It is easy
to see that $\tilde m_\pi(1^i) = i^{\underline{\ell(\pi)}}$,
whence $\X_G(1^i) = \chi_G(i)$, so that $\X_G$ is a generalization
of the chromatic polynomial.
Stanley proves a number of results about~$\X_G$ in the papers [St2][St4]
(see also [Ga1][Ga2]).
We shall be investigating some other properties of~$\X_G$ later.

Now let $D$ be a digraph.  Following Chung and Graham~[CG1],
we say that a subset~$S$ of the edges
of~$D$ is a \itc{path-cycle cover} of~$D$ if no two elements of~$S$ lie
in the same row or column of the associated
board.\note{Chung and Graham have a revised version [CG2] of the
preprint [CG1].  We shall be referring to both versions.}
If we think of~$S$ as a spanning
subgraph of~$D$ then we see that this condition just means that $S$ is
a (vertex-)disjoint union of directed paths and directed cycles.
(Isolated vertices are thought of as directed paths with zero edges,
and isolated loops are thought of as cycles of length one.)
A path-cycle cover with no cycles is called a \itc{path cover,}
and a path-cycle cover with no paths is called a \itc{cycle cover.}
The \itc{type} of a path-cycle cover~$S$ is the
set partition of~$V(D)$ such that each block is the set of vertices
of one of these directed paths or directed cycles.  We write $\pi(S)$
for the set of blocks corresponding to the directed paths and
$\sigma(S)$ for the set of blocks corresponding to the directed cycles,
and we say that the type of~$S$ is $(\pi,\sigma)$ if $\pi(S)=\pi$ and
$\sigma(S)=\sigma$.
Chung and Graham's \itc{cover polynomial~$C(D;i,j)$} is then defined by
$$C(D; i,j) \defeq
   \sum_S i^{\underline{\ell(\pi(S))}} j^{\ell(\sigma(S))},$$
where the sum is over all path-cycle covers $S\subset E(D)$.

In view of these definitions and the fact that $p_\sigma(1^j) =
j^{\ell(\sigma)}$, the following definition (suggested by Stanley~[St1])
is quite natural.

\definition
Let $D$ be a digraph, and let
$x = \{x_1, x_2, \ldots\}$ and $y=\{ y_1, y_2, \ldots\}$
be two sets of commuting independent indeterminates.
The \itc{path-cycle symmetric function~$\Xi_D$} of~$D$ is defined by
$$\Xi_D = \Xi_D(x;y)
   \defeq \sum_S \tilde m_{\pi(S)}(x) \, p_{\sigma(S)}(y),$$
where the sum is over all path-cycle covers $S\subset E(D)$.

The path-cycle symmetric function has not been investigated before, and
the bulk of this thesis is devoted to its study.

Note that if we only care about \itc{path} covers we can simply
consider $\Xi_D(x;0)$.  In addition, if $B$ is the board associated
with~$D$, then $\Xi_D(0;y)$ is equivalent to what Stanley and Stembridge
call~$Z[B]$ ([S-S, section~3]).  Thus
we may regard $\Xi_D$ as a further generalization of Stanley and
Stembridge's generalization of the theory of permutations with
restricted position.

The following fact is immediate.

\proposition\restrict{$\Xi_D(1^i;1^j) = C(D;i,j)$.\qed}

There is a close connection between $\X_G$ and~$\Xi_D$.
Given a poset~$P$, let $G(P)$ denote its incomparability graph (in which
two vertices of the poset are adjacent if and only if
they are incomparable), and let
$D(P)$ denote the digraph with edge set $\{(i,j) \mid i<j\}$.  Chung and
Graham observe that for any poset~$P$,
$$C\bigl(D(P); i, 0\bigr) = \chi_{G(P)}(i).$$
This connection generalizes readily to the symmetric function case.

\proposition\connection{For any poset~$P$, $\Xi_{D(P)} = \X_{G(P)}$.}

\proof
Since $D(P)$ is acyclic, all path-cycle covers are in fact just path
covers, so the $y$~variables can be deleted from the definition of~$\Xi_D$
in this case.  But path covers of~$D(P)$ correspond to partitions of~$P$
into chains, which correspond to stable partitions of~$G(P)$.  Comparing
the definitions of $\Xi_D$ and~$\X_G$ yields the proposition.\qed

\sec Connection with rook theory

Let $B\subset [d]\times[d]$ be a board, and let the
\itc{rook number~$r_k^B$} denote the number of
ways of placing $k$ non-taking rooks on~$B$ (i.e., the number of subsets
of~$B$ such that no two squares lie in the same row or in the same column).
Following Goldman, Joichi and
White~[GJW1], we define the \itc{$d$-factorial polynomial}
(or simply the \itc{factorial polynomial}) $R(B;i)$ by
$$R(B;i) \defeq \sum_k r_k^B i^{\underline{d-k}}.$$
If $D$ is the digraph associated with~$B$, we also write $r_k^D$ for~$r_k^B$
and $R(D;i)$ for~$R(B;i)$.
(With this equivalence between boards and digraphs, the
factorial polynomial is the same as Chung and Graham's
\itc{binomial drop polynomial.})
The study of the factorial polynomial and other rook polynomials
is a well-established area of combinatorics (see for example
[Rio, Chapters 7 and~8][GJW1][GJRW][GJW2][GJW3][GJW4]).
The definition of a path-cycle cover already suggests a connection with rook
theory.  More precisely, we have the following proposition.

\proposition\coverrook{For any digraph~$D$,
$R(D;i)=C(D;i,1)=\Xi_D(1^i;1)$.}

\proof
The first equality is demonstrated in~[CG1] (or~[CG2]) and the second
equality follows from \restrict.\qed

\coverrook\
(as well as, for example, [S-S, section~3] and [St2, Proposition~5.5])
suggests that some of the theory of rook polynomials might generalize
to~$\Xi_D$.  This is indeed the case, as we shall see in more detail
later.


\chap The Path-Cycle Symmetric Function

In this chapter, we carry out a fairly systematic investigation of
the path-cycle symmetric function.  For example, we try to derive
as many analogues of theorems about~$\X_G$ and generalizations of
results from rook theory as we can.  In addition, as we shall see,
we uncover some unexpected results as a byproduct of this investigation.

\sec Basic facts

It is always good practice to begin the study of any
mathematical object by computing a few examples.
The chromatic symmetric functions of some graphs are computed in~[St2],
and \connection\ tells us that when these graphs are incomparability
graphs, the chromatic symmetric functions are also path-cycle symmetric
functions.  But it would also be nice to see some examples of path-cycle
symmetric functions that do not arise in this way.

In general, computing the path-cycle symmetric function is very hard
(in~the sense of
computational complexity).  Even the computation of the cover polynomial
is \#$P$-hard, since, as Chung and Graham observe, $C(D; 1,0)$ is the
number of Hamiltonian paths, and computing this number is
\#$P$-hard (an~indication of how this last fact may be proved can be found
in section~7.3 of~[G-J]).
Furthermore, there are not many digraphs~$D$ for which
$\Xi_D$ has a simple closed form.
For example, since $\Xi_D$ is
a generalization of Goldman, Joichi and White's factorial polynomial,
for which there is a nice factorization theorem
in the case of Ferrers shapes~([GJW1]),
we might hope for a factorization theorem for~$\Xi_D$.
Unfortunately, this is not the case;
the factorization theorem does not even generalize straightforwardly
to the cover polynomial (although Dworkin~[Dwo] has found some cases
where the cover polynomial does factorize).
Part of the problem is that unlike the factorial polynomial, $\Xi_D$
(or even the cover polynomial) does not remain invariant under permutations
of rows and columns.
Furthermore, like $\X_G$, the path-cycle symmetric function does not satisfy
a deletion-contraction recurrence.

There are, however, some nice formulas for the path-cycle symmetric
functions of directed paths and directed cycles.

\proposition\dirpathcycle{Let $P_d$ be the directed path with $d$ vertices
and let $C_d$ be the directed cycle with $d$ vertices.  Then
$$\eqalign{
 \Xi_{P_d} &= \sum_{\lambda\vdash d} \ell(\lambda)!\, m_\lambda(x)
  = \sum_{r=0}^{d-1} u_r s_{d-r,1^r}(x)\cr
{\rm and}\qquad
 \Xi_{C_d} &= d\sum_{\lambda\vdash d} \bigl(\ell(\lambda) - 1\bigr)!\,
  m_\lambda(x) + p_d(y) = d\sum_{r=0}^{d-1} v_r s_{d-r,1^r}(x) +
  p_d(y),\cr}$$
where $u_r$ is the number of permutations of $r+1$ with no consecutive
ascending pairs, and $v_r$ is the number of permutations of~$r$ with
no fixed points.}

\noindent
(Remark: the Schur function expansions were obtained with the aid of~[Slo].)

\proof
Observe that \itc{every} subset~$S$ of the
edges of~$P_d$ is a path cover.  Now $\pi(S)$ has type~$\lambda$ if and
only if the complement of~$S$ is a set of edges whose removal from~$P_d$
breaks $P_d$ into paths whose sizes are given by
the parts of~$\lambda$.  So the number of such sets~$S$ is
$${\ell(\lambda) \choose r_1, r_2, \ldots},$$
where $r_i$ is the number of parts of~$\lambda$ of size~$i$.  Hence
$$\Xi_{P_d} = \sum_{\lambda\vdash d} {\ell(\lambda) \choose r_1, r_2, \ldots}
  \tilde m_\lambda = \sum_{\lambda \vdash d} \ell(\lambda)!\, m_\lambda.$$
To establish the Schur function expansion of~$\Xi_{P_d}$,
first note that
$$\sum_{r=0}^{k-1} u_r {k-1 \choose r} = k!$$
(see [Kre] or~[R-P] for a proof).
Thus it suffices to show that the coefficient of
$m_\lambda$ in $s_{d-r,1^r}$ is
$$\ell(\lambda) - 1 \choose r$$
whenever $\lambda \vdash d$, for then it will follow that
$$\sum_{r=0}^{d-1} u_r s_{d-r,1^r}
  = \sum_{\lambda\vdash d} m_\lambda
     \sum_{r=0}^{d-1} u_r {\ell(\lambda) - 1 \choose r}
  = \sum_{\lambda\vdash d} \ell(\lambda)!\, m_\lambda,$$
since $d\ge\ell(\lambda)$ for all $\lambda\vdash d$.

The coefficient of $m_\lambda$ in $s_{d-r,1^r}$ is just the number
of column-strict Young tableaux of type $(d-r,1^r)$ and content~$\lambda$
(see~[Mac, (5.12)]).  These Young tableaux are precisely those obtained by
the following procedure: put a~1 in the corner of the tableau, and then
choose any $r$ distinct integers in the set $\{2,3,\ldots,\ell(\lambda)\}$
and put them in order down the left column of the tableau.
Then fill in the rest of the first row with whatever integers are needed
to make the content equal to~$\lambda$ (this can be done in exactly one
way since the row must increase monotonically from left to right).
Hence the number of such tableaux is
$${\ell(\lambda) - 1 \choose r},$$
as required.

To compute $\Xi_{C_d}$, observe that the set of \itc{all} edges of~$C_d$
is a cycle cover (this accounts for the term~$p_d(y)$) and that every other
subset of edges is a path cover.  The following procedure generates
all path covers of~$C_d$ of type~$\lambda$:
take a set of directed paths whose sizes are given by the parts of~$\lambda$,
and make paths of the same length distinguishable; then arrange the paths to
form a circle, and choose some vertex to be vertex~1.
Now there are $\bigl(\ell(\lambda)-1\bigr)!$ distinct
circular permutations of the paths, and $d$ ways to choose a vertex~1,
but notice that this procedure generates each path cover
$r_1!r_2!\cdots$ times because ``in reality''
directed paths of the same length are not distinguishable.
This gives us the first formula for~$\Xi_{C_d}$.
To prove the other formula, all we need to establish (in view of the
remarks above in the case of~$P_d$) is that
$$\sum_{r=0}^{k-1} v_r {k - 1 \choose r} = (k-1)!,$$
but this is easy: an arbitrary permutation of $k-1$ letters can be chosen
by first choosing $k-1-r$ fixed points and then choosing a permutation of
the remaining $r$~letters that has no fixed points.\qed

We have found two formulas which allow the path-cycle symmetric
function of a digraph to be computed from path-cycle symmetric functions
of related digraphs.  One of these formulas is sufficiently interesting
that we devote the entire next section to it.  The other formula is a
multiplicativity property that generalizes Chung and Graham's
multiplicativity formula for the cover polynomial [CG2, Corollary~2].
To prove this result,
we introduce the concept of a path-cycle coloring,
due to Chung and Graham (see [CG2];
the definition in~[CG1] contains a minor error).

\definition
A \itc{path-cycle coloring} of a digraph~$D$ is an ordered
pair $(S,\kappa)$ where $S$ is a path-cycle cover and $\kappa$ is a
map from~$V(D)$ to the positive integers such that

\item{1.} $\kappa(v_1) = \kappa(v_2)$ if $v_1$ and~$v_2$ belong to the same
path or if $v_1$ and~$v_2$ belong to the same cycle, and
\item{2.} $\kappa(v_1) \ne \kappa(v_2)$ if $v_1$ belongs to a path
and $v_2$ belongs to a \itc{different} path.

A \itc{path coloring} is a path-cycle coloring with no cycles.

\proposition\coloring{For any digraph~$D$,
$$\Xi_D = \sum_{(S,\kappa)} \;
  \prod_{u\rm\;is\;in\;a\;path} \!\!\! x_{\kappa(u)}
  \prod_{v\rm\;is\;in\;a\;cycle} \!\!\! y_{\kappa(v)},$$
where the sum is over all path-cycle colorings~$(S,\kappa)$.}

\proof
Regard the sum as a double sum: for each path-cycle cover~$S$, sum over
all ``compatible'' colorings~$\kappa$, and then sum over all~$S$.  For
each fixed~$S$ the paths and cycles may be colored independently so the
sum over~$\kappa$ factors into a product of a symmetric function in~$x$
and a symmetric function in~$y$.  Clearly coloring the paths with
distinct colors gives~$\tilde m_{\pi(S)}(x)$ and coloring the cycles
so that each cycle is monochromatic gives~$p_{\sigma(S)}(y)$.\qed

\proposition\multiplicativity{Suppose $D$ is the digraph formed by joining
the disjoint digraphs $D_1$ and~$D_2$ with all the edges $(v_1,v_2)$
with $v_1\in V(D_1)$ and $v_2\in V(D_2)$.  Then
$\Xi_D = \Xi_{D_1}\Xi_{D_2}$.}

\proof
A path-cycle coloring of~$D$ induces path-cycle colorings of both
$D_1$ and~$D_2$ by restriction.
Conversely, given any path-cycle coloring of~$D_1$
and any path-cycle coloring of~$D_2$ there exists a unique path-cycle
coloring of~$D$ inducing them: if a path in~$D_1$ has the same color
as a path in~$D_2$, join them end to end with the appropriate edge from
$D_1$ to~$D_2$.  The result now follows from \coloring.\qed

Notice that if $D_1=D(P_1)$ and~$D_2=D(P_2)$ for some posets $P_1$
and~$P_2$, then the construction described in \multiplicativity\
corresponds to the \itc{ordinal sum} $P_1\oplus P_2$
(see~[St3, Chapter~3]).
Thus \multiplicativity\ can save us some labor in computing $\Xi_{D(P)}$
if $P$ has a nontrivial ordinal sum decomposition.  For example, it is
clear from the definitions that if $D$ has no edges at all, then
$\Xi_D=\tilde m_{1^d} = d!e_d$.  It follows that
$\Xi_{D(P)}$ is a multiple of an elementary symmetric function
whenever $P$ is an ordinal sum of antichains.

We now turn from the problem of computing~$\Xi_D$ to the problem of
finding interesting facts about it.  A natural thing to do is to try
and find analogues of known theorems about~$\X_G$.  This strategy is
not always successful, but we do have the following
analogues of Corollaries 2.7 and~2.11 of~[St2].

\proposition\ppositivity{If $D$ is an acyclic digraph, then
$\omega_x\Xi_D$ is $p$-positive.}

\proof
Since $D$ is acyclic, all path-cycle covers are path covers,
and $\Xi_D = \Xi_D(x,0)$.  From
Doubilet~[Dou, Appendix~1] we know that for any set partition~$\pi$,
$$\tilde m_\pi = \sum_{\sigma\ge\pi} \mu(\pi,\sigma) p_\sigma.$$
Thus
$$\eqalign{\Xi_D &= \sum_S \sum_{\sigma\ge\pi(S)}
      \mu\bigl(\pi(S),\sigma\bigr) p_\sigma\cr
  &= \sum_\sigma \biggl(\sum_{\{S\mid \pi(S)\le\sigma\}}
      \mu\bigl(\pi(S),\sigma\bigr) \biggr) p_\sigma,\cr}$$
where $S$ ranges over path covers.
Now fix $\sigma$ and let $D_1, D_2,\ldots, D_l$ be the subgraphs induced by
the blocks of~$\sigma$, with sizes $d_1,d_2,\ldots,d_l$ respectively.
If $c_i$ is the coefficient of $p_{d_i}$ in~$\Xi_{D_i}$, then
we claim that the coefficient of~$p_\sigma$ in~$\Xi_D$ is $\prod_i c_i$.
To see this, first note that choosing a path cover~$S$ of~$D$ such 
that $\pi(S)\le\sigma$ is equivalent to (independently) choosing path
covers for each~$D_i$.
If we let $\Pi_n$ denote the lattice of partitions of~$[n]$ ordered
by refinement, then it is well known and easy to see that the interval $[\hat 0,\sigma]$ in a partition lattice
is isomorphic to
$$\Pi_{d_1} \times \Pi_{d_2} \times \cdots \times \Pi_{d_l}.$$
It is also well known (e.g.,~[St3, Prop.~3.8.2])
that the M\"obius function of a product
is the product of the M\"obius functions.
Putting these facts together readily yields our claim.

Thus to prove the theorem it suffices to prove that for an acyclic digraph
with $d$~vertices the sign of the coefficient of~$p_d$ is~$(-1)^{d-1}$.
For then, since any induced subgraph of an acyclic graph is acyclic, we
can apply our claim above to show that the coefficient of~$p_\sigma$
is~$\sgn\sigma$.

Let $D$ have $d$ vertices.
By specializing via \coverrook, we see that the
coefficient of~$p_d$ in $\Xi_D$ equals the coefficient of~$i$ in $R(D;i)$.
Directly from the definitions we see that this equals
$$(-1)^{d-1}\sum_{k=0}^{d-1} (-1)^k r_k^D (d-k-1)!$$
Now an acyclic digraph has at least one source and one sink, so by removing
the corresponding row and column we see that $B$ may be regarded as a
subset of a $(d-1)\times(d-1)$ board.  Then the above sum is a positive
integer, by the inclusion-exclusion formula for rooks
(see~[St3, Theorem~2.3.1]).
The sign of the coefficient is therefore~$(-1)^{d-1}$ as desired.\qed

We mention in passing that the argument used in the above proof allows
[St2, Proposition~5.5] to be extended from posets to acyclic digraphs.

\proposition\partialderiv{If $D$ is an acyclic digraph, then
$${\partial\Xi_D\over\partial p_i} = \sum_H
  \biggl(\sum_S \mu\bigl(\pi(S),V(H)\bigr)\biggr)\Xi_{D\backslash H},$$
where $H$ runs over all (vertex-)induced subgraphs of~$D$ with $i$ vertices,
$S$ runs over all path covers of~$H$, $\mu$ is the M\"obius function of the
lattice of partitions of~$V(H)$, and $D\backslash H$ is the subgraph
obtained by deleting $V(H)$ from~$D$.}

\proof
For the purposes of this proof only, we use the notation $S\sqsubseteq H$
to mean ``$S$~is a path cover of~$H$.''
From the proof of \ppositivity, we have
$$\Xi_D = \sum_\sigma
  \biggl(\sum_{\{S\sqsubseteq D\mid\pi(S)\le\sigma\}}
  \mu\bigl(\pi(S),\sigma\bigr) \biggr) p_\sigma,$$
where $\sigma$ ranges over all partitions of~$V(D)$.
Differentiating both sides with respect to~$p_i$, we obtain
$${\partial\Xi_D\over\partial p_i} = \sum_{(H,\tau)}
  \biggl(\sum_{\bigl\{S\sqsubseteq D\bigm|\pi(S)\le\tau\cup\{V(H)\}\bigr\}}
  \mu\bigl(\pi(S),\tau\cup\{V(H)\}\bigr) \biggr) p_\tau,$$
where the outer sum ranges over all ordered pairs $(H,\tau)$ such that
$H$ is an induced subgraph of~$D$ with $i$ vertices and $\tau$ is any
partition of $V(D\backslash H)$.  Since, as noted in the proof of
\ppositivity, the M\"obius function is multiplicative, it follows that
$$\eqalign{{\partial\Xi_D\over\partial p_i} &= \sum_{(H,\tau)}
  \biggl(\sum_{S\sqsubseteq H} \mu\bigl(\pi(S),V(H)\bigr)\biggr)
  \biggl(\sum_{\{S \sqsubseteq D\backslash H \mid \pi(S) \le \tau\}}
     \mu\bigl(\pi(S),\tau\bigr)\biggr) p_\tau \cr
 &= \sum_H
  \biggl(\sum_{S\sqsubseteq H} \mu\bigl(\pi(S),V(H)\bigr)\biggr)
    \sum_\tau
  \biggl(\sum_{\{S \sqsubseteq D\backslash H \mid \pi(S) \le \tau\}}
     \mu\bigl(\pi(S),\tau\bigr)\biggr) p_\tau \cr
 &= \sum_H
  \biggl(\sum_{S\sqsubseteq H} \mu\bigl(\pi(S),V(H)\bigr)\biggr)
    \Xi_{D\backslash H}.\qed\cr}$$

We conclude this section with two counterexamples.  The digraphs
\vskip 0.3truein
\hskip 0.6truein
\psfig{figure=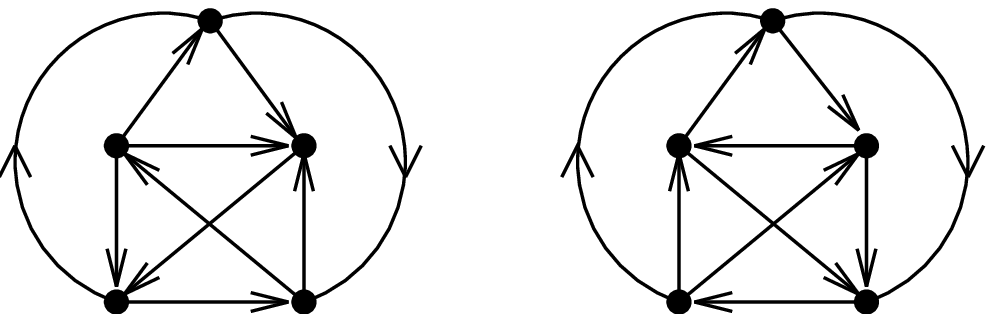}

\noindent
(taken from figure~2 of~[Bon])
have identical lists of vertex-deleted subdigraphs but have
different path-cycle symmetric functions (and even different
cover polynomials), and hence the path-cycle symmetric function
is not reconstructible (as opposed to the chromatic symmetric
function---see section~3 of the next chapter).

The reader who is aware that
every factorial polynomial is the chromatic polynomial
of some graph~([GJW2])
might wonder if this generalizes to our symmetric function context.
The answer is no, and
the directed path on three vertices provides an example of
an acyclic digraph~$D$ for which $\Xi_D$ does not equal $\X_G$ for any
graph~$G$.
The philosophical reason for this is that the proof that every factorial
polynomial is a chromatic polynomial
relies on the fact (mentioned previously)
that the factorial polynomial of a board is
independent of how the board is embedded in its $[d]\times[d]$ grid,
but the same cannot be said of the path-cycle symmetric function.

\sec Reciprocity

The \itc{complement~$D'$} of a digraph~$D$ is the digraph on
the same vertex set whose edges are precisely those pairs $(i,j)$ that
are not edges of~$D$.  In this section we prove one of the most striking
facts about the path-cycle symmetric function; namely, a combinatorial
reciprocity theorem relating
$\Xi_D$ and~$\Xi_{D'}$.  We shall need two
change-of-basis formulas, which we shall now state.

As in the proof of \ppositivity,
let $\Pi_n$ denote the lattice of partitions of~$[n]$
(ordered by refinement).  Recall that if $\pi\le\sigma$ in~$\Pi_n$
and $r_i$ is the number of blocks of~$\sigma$ that are composed of
$i$~blocks of~$\pi$, then the M\"obius function satisfies
$$|\mu(\pi,\sigma)| = \prod_i (i-1)!^{r_i}.$$
(See [St3, Example~3.10.4] for a proof.)
Also, following Doubilet~[Dou], define
$$\lambda(\pi,\sigma)! \defeq \prod_i i!^{r_i}.$$
We then have the following change-of-basis formulas (taken from~[Dou,
Appendix~1]).

\proposition\doubilet{
$$f_\pi = \sum_{\sigma\ge\pi} \lambda(\pi,\sigma)! \, \tilde m_\sigma
        = \sum_{\sigma\ge\pi} |\mu(\pi,\sigma)| \, p_\sigma.\qed$$}

We are now ready for the main theorem of this section.

\theorem\reciprocity{For any digraph $D$,
$$\Xi_D(x;y) = \sum_S \sgn\pi(S) f_{\pi(S)}(x,y)\,
    p_{\sigma(S)}(-y),$$
where the sum is over all path-cycle covers of the complement~$D'$.
Equivalently,
$$\Xi_D(x;y) = [\omega_x \Xi_{D'}(x;-y)]_{x\to(x,y)},$$
where $[g(x;y)]_{x\to(x,y)}$ means that, treating $g$ as a symmetric
function in the~$x$'s with coefficients in the~$y$'s,
the set of $x$~variables
is to be replaced by the union of the $x$ and~$y$ variables.}

\noindent
(Remark: ``replacing the $x$ variables with the union of the $x$ and~$y$
variables'' may be formalized as ``applying the endomorphism~$\Delta$ of
$\Lambda(x)\otimes\Lambda(y)$ that sends $p_m(x)$ to $p_m(x)+p_m(y)$ and
that leaves $p_n(y)$ fixed.'')

\proof
The equivalence of the two formulations is clear.
We define a \itc{partitioned order} of~$D$ to be a partition of~$V(D)$
together with either a linear order or a cyclic order on each block.
If $\kappa$ is a partitioned order of~$D$, let $\pi(\kappa)$ be the
set of blocks with linear orders and let $\sigma(\kappa)$ be
the set of blocks with cyclic orders.  Let $E_\kappa$ denote
the set of ordered pairs~$(u,v)$ satisfying the following two conditions.

\item{1.} $u$ and~$v$ are in the same block of~$\kappa$ and $u$
immediately precedes~$v$ in the linear or cyclic order on the block.

\item{2.} $(u,v)$ is not an edge of~$D$.

Note that $E_\kappa \subset E(D')$ and that there is a natural bijection
between partitioned orders $\kappa$ such that $E_\kappa = \varemptyset$
and path-cycle covers
(given such a partitioned order, take all $(u,v)$
satisfying condition~1 above).
Now for any finite set~$T$, the alternating sum
$$\sum_{S \subset T} (-1)^{|S|}$$
equals one if $T = \varemptyset$ and is zero otherwise.  Thus
$$\Xi_D = \sum_\kappa \tilde m_{\pi(\kappa)}\,p_{\sigma(\kappa)}
          \sum_{S \subset E_\kappa} (-1)^{|S|},$$
where the first sum is over all partitioned orders of~$D$.  We now
interchange the order of summation.  Observe first that all sets~$S$
that arise are path-cycle covers of~$D'$, since $S$ is a subset of
the set of all $(u,v)$ satisfying condition~1 above for some~$\kappa$.
Given a path-cycle cover~$S$
of~$D'$, we now need to determine the set~$\Scr P$
of partitioned orders of~$D$ that give rise to it.
Only blocks with cyclic orders can give rise to cycles of~$S$,
so for every $\kappa\in\Scr P$, $\sigma(\kappa)$ must include the
blocks of~$\sigma(S)$ among its own blocks.  On the other hand, the
blocks of~$\pi(S)$ can arise either from blocks with linear orders
or from blocks with cyclic orders.  To determine all possibilities
we must consider all ways of agglomerating the blocks of~$\pi(S)$
into blocks of~$\pi(\kappa)$, and then for each composite block
in each such agglomeration we must consider both linear and cyclic
orders.  The linear or cyclic order on the composite block can be
viewed as a linear or cyclic order on the \itc{blocks of~$\pi(S)$}
(instead of on the vertices), because the linear or cyclic order must
induce the edges of~$S$, i.e., if $(u,v)$ is an edge of~$S$ then $u$ must
immediately precede~$v$ in the order dictated by~$\kappa$, and therefore
the vertices in each block of~$\pi(S)$ are constrained to be consecutive
and in a fixed order.
Clearly, every such linear or cyclic order
on the blocks gives rise to a unique $\kappa\in\Scr P$.  The number of
ways to
impose a linear order if there are $i$~blocks is $i!$ and the number of
ways to impose a cyclic order is $(i-1)!$.  Thus we can enumerate~$\Scr P$
by summing over all divisions of the blocks of~$\pi(S)$ into two groups
$\alpha$ and~$\beta$
(linear and cyclic) and, for each such division, summing over all ways
of grouping the blocks into composite blocks, weighted by a factorial
factor.  More precisely we have
$$\Xi_D = \sum_S (-1)^{|S|} p_{\sigma(S)}(y) \sum_{(\alpha,\beta)}
 \sum_{\scriptstyle\gamma\ge\alpha\atop\scriptstyle\delta\ge\beta}
    \lambda(\alpha,\gamma)!\,|\mu(\beta,\delta)|\,\tilde m_{\gamma}(x)\,
    p_{\delta}(y),$$
where the first sum is over all path-cycle covers of~$D'$.
By \doubilet, we have
$$\eqalign{\Xi_D &= \sum_S (-1)^{|S|} p_{\sigma(S)}(y) \sum_{(\alpha,\beta)}
   f_{\alpha}(x)\, f_{\beta}(y) \cr
   &= \sum_S (-1)^{|S|} p_{\sigma(S)}(y) \sum_{(\alpha,\beta)}
   (\sgn\alpha)(\sgn\beta)\, \omega_x \tilde m_{\alpha}(x)\,
   \omega_y \tilde m_{\beta}(y).\cr}$$
Now the blocks of $\alpha$ and~$\beta$ correspond to the paths of~$S$,
so $|\alpha| - \ell(\alpha)$ is the number of edges of~$S$ in~$\alpha$,
and similarly for~$\beta$.  Thus $(\sgn\alpha)(\sgn\beta)$ depends only
on the total number of edges of~$S$ devoted to directed paths
(namely, $|\pi(S)|-\ell(\pi(S))$) and does
not depend on the particular choice of~$\alpha$ or~$\beta$.  We have
$$\eqalign{\Xi_D &= \sum_S (-1)^{|S|} p_{\sigma(S)}(y)
   \, (-1)^{|\pi(S)|-\ell(\pi(S))} \sum_{(\alpha,\beta)}
   \omega_x \tilde m_{\alpha}(x)\,
   \omega_y \tilde m_{\beta}(y)\cr
   &= \sum_S (-1)^{|\sigma(S)|} p_{\sigma(S)}(y) \,\omega_x \omega_y\!
   \sum_{(\alpha,\beta)} \tilde m_{\alpha}(x)\, \tilde m_{\beta}(y).\cr}$$
A moment's thought shows that the inner sum is $\tilde m_{\pi(S)}(x,y)$.
Now
$$\eqalign{\omega_x\omega_y p_n(x,y) &=
  \omega_x\omega_y \bigl(p_n(x)+p_n(y)\bigr)
 = (-1)^{n-1}p_n(x) + (-1)^{n-1} p_n(y) \cr
 & = [\omega_x p_n(x)]_{x\to(x,y)},\cr}$$
and since $\omega$ and $x\to(x,y)$ are both homomorphisms,
$$\omega_x\omega_y g(x,y) = [\omega_x g(x)]_{x\to(x,y)}$$
for any symmetric function~$g$.
Finally, $(-1)^{|\sigma|}p_{\sigma}(y) = p_{\sigma}(-y)$, so we obtain
$$\Xi_D = \sum_S p_{\sigma(S)}(-y)\,
   [\omega_x \tilde m_{\pi(S)}(x)]_{x\to(x,y)}$$
as desired.\qed

We remark that the appearance of $\omega$ in \reciprocity\
is what leads us to call it a combinatorial reciprocity theorem.

\reciprocity\ is a rather curious result in that it is not obvious
that the operation
$$\iota: g(x;y) \to [\omega_x g(x;-y)]_{x\to(x,y)}$$
is an involution.  One might even wonder if $\iota$ is really
an involution on \itc{all} of $\Lambda(x)\otimes\Lambda(y)$
or if it is only an involution when restricted to some subset
of $\Lambda(x)\otimes\Lambda(y)$
that includes all path-cycle symmetric functions.  It is easily
verified, however, that
$$\iota\bigl(\iota\bigl( p_m(x) p_n(y) \bigr)\bigr) = p_m(x)p_n(y)$$
for any $m$ and~$n$.
Furthermore, $\iota$ is the composition of three
operations, each of which is an endomorphism of
$\Lambda(x)\otimes\Lambda(y)$.
It follows that $\iota$ is an involution on all of
$\Lambda(x)\otimes\Lambda(y)$.

On the other hand, one can still ask whether there is any way of
reformulating \reciprocity\ in a way that makes it more obvious
that the operation involved is an involution.  One way of doing
this is to define
$$\hat \Xi_D(x;y) \defeq \sum_S (-2)^{\ell(\sigma(S))}
    \tilde m_{\pi(S)}(x,y)\, p_{\sigma(S)}(y),$$
where the sum is still over all path-cycle covers $S$ of~$D$.
Then we have the following result.

\proposition\reciprocityhat{$\hat \Xi_{D'}(x;y) =
   \omega_x \hat \Xi_D(x;-y)$.}

\noindent
(Remark: I arrived at the mysterious-looking
definition of $\hat\Xi_D$ by first
specializing to the cover polynomial, finding a change of
variables that made the reciprocity operation an obvious
involution, and then generalizing back to the symmetric
function case by using plethysm.)

\proof
Let $\kappa$ be the endomorphism of $\Lambda(x)\otimes\Lambda(y)$
that sends $p_n(y)$ to $-2p_n(y)$ and that leaves $p_m(x)$ fixed.
Then $\hat\Xi_D$ is obtained from $\Xi_D$ by applying $\kappa$ first
and then $\Delta$.
In view of \reciprocity, it suffices to prove that the
following operations have the same effect on all elements of
$\Lambda(x)\otimes\Lambda(y)$:

\item{1.} (Left-hand side.)  Apply $\iota$, then $\kappa$, then $\Delta$.

\item{2.} (Right-hand side.) Apply $\kappa$, then apply $\Delta$, then
negate the $y$ variables, and then apply~$\omega_x$.

Each of these operations is an endomorphism, so it suffices to check their
behavior on $p_m(x)p_n(y)$.  A straightforward computation shows that in
both cases the result is
$$2p_n(y)\bigl(p_m(x) - p_m(y)\bigr)(-1)^{m+n},$$
which proves the desired result.\qed

While \reciprocityhat\ has the advantage of involving an operation that
is clearly an involution, it has the disadvantage that the definition
of $\hat\Xi_D$ is mysterious.  In particular, the $-2$ has no obvious
combinatorial significance.  For this reason, we regard $\hat\Xi_D$
as an artificial contrivance and we shall continue to use $\Xi_D$
instead.

\reciprocity\ readily yields several attractive corollaries.
For example, by setting all the $x$ variables equal to zero, we
immediately obtain [S-S, Theorem~3.2].  More interestingly, we can
obtain an affirmative answer to the question, raised by Chung and
Graham [CG1, section 8(c)], of whether
$C(D; i,j)$ determines $C(D'; i,j)$.

\corollary\coverreciprocity{If $D$ is a digraph with $d$ vertices, then
$C(D'; i,j) = (-1)^d C(D;{}-i-j,j)$.}

\proof
Let $g$ be any symmetric function that is homogeneous of degree~$d$
and let $g^* = \omega g$.
We claim that $g^*(1^i)$ is obtained by changing $i$ to~$-i$ in~$g(1^i)$
and then multiplying by~$(-1)^d$.  To see this, first consider the case
where $g=p_\lambda$ for some $\lambda\vdash d$.  Then
$g^* = (\sgn\lambda)p_\lambda$ and hence
$$g^*(1^i) = (\sgn\lambda) i^{\ell(\lambda)}.$$
On the other hand $g(1^i) = i^{\ell(\lambda)}$.  Changing $i$
to~$-i$ and multiplying by~$(-1)^d$ amounts to multiplying
by~$(-1)^{d-\ell(\lambda)} = \sgn\lambda$, as required.  The claim then
follows by linearity.

Now $\tilde m_\pi(1^i,1^j) = (i+j)^{\underline{\ell(\pi)}}$.  Since
$(\sgn\pi) f_\pi = \omega\tilde m_\pi$, we have
$$(\sgn\pi) f_\pi(1^i, 1^j) = (-1)^{|\pi|} ({}-i-j)^{\underline{\ell(\pi)}}.$$
Also, as noted before, $p_\sigma(-y) = (-1)^{|\sigma|} p_\sigma(y)$.  Thus,
specializing \reciprocity\ via \restrict\ yields
$$\eqalign{C(D; i,j) &= \sum_S (-1)^{|\sigma(S)|} (-1)^{|\pi(S)|}
            ({}-i-j)^{\underline{\ell(\pi(S))}} j^{\ell(\sigma(S))}\cr
    &= (-1)^d C(D'; {}-i-j,j).\qed\cr}$$

\coverreciprocity\ can be proved directly using
deletion-contraction techniques, and it has also been obtained
independently by Gessel~[Ge1].  We omit the details.

A further specialization of \reciprocity\ gives a formula for rook
polynomials; we defer this to the next section, where we consider rook
theory in more detail.

\corollary\onevarreciprocity{For any digraph~$D$,
$$\Xi_D(x;0) = \omega_x \Xi_{D'}(x;0).$$}

Note the similarity between this result and Stanley's reciprocity
theorem [St2, Theorem~4.2].  In fact, the two reciprocity theorems overlap,
because of \connection, so \onevarreciprocity\ gives a new interpretation
of $\omega \Xi_{D(P)} = \omega \X_{G(P)}$ when $P$ is a poset.

\onevarreciprocity\
follows immediately from \reciprocity, but we shall give two other
proofs because they illustrate connections with other known results.
The first proof is due to Gessel~[Ge1], and it derives
\onevarreciprocity\ from
a result of Carlitz, Scoville and Vaughan [CSV, Theorem~7.3].
We need some preliminaries.  Given a digraph~$D$ with $d$ vertices, let
$$A_D = \{a_1,a_2,\ldots,a_d\}$$
be a set of commuting independent indeterminates, and define
$$\alpha_{D,n} = \sum_{i_1,i_2,\ldots,i_n} a_{i_1}a_{i_2}\cdots a_{i_n},$$
where the sum is over all $i_1,i_2,\ldots,i_n$ such that
$(a_{i_j},a_{i_{j+1}})$ is an edge of~$D$ for all $j<n$.  Similarly, let
$$\alpha'_{D,n} = \sum_{i_1,i_2,\ldots,i_n} a_{i_1}a_{i_2}\cdots a_{i_n},$$
where this time the sum is over all $i_1,i_2,\ldots,i_n$ such that
$(a_{i_j},a_{i_{j+1}})$ is an edge of the complement~$D'$ for all $j<n$.
With this notation, the result of Carlitz, Scoville and Vaughan is
(essentially) the following.

\proposition\carlitz{For any digraph~$D$,
$$\sum_n (-1)^n \alpha'_{D,n} =
  \biggl(\sum_n \alpha_{D,n} \biggr)^{-1}.\qed$$}

We can now give Gessel's proof of \onevarreciprocity.

\noindent{\it First proof of \onevarreciprocity.}
Let $\theta_{D,y}$ be the homomorphism from the ring of symmetric functions
in the variables $y=\{y_1,y_2,\ldots\}$ to the ring of formal power series
in~$A_D$ that sends the complete symmetric function $h_n(y)$
to~$\alpha_{D,n}$.  Similarly, let $\theta'_{D,y}$ be the homomorphism that
sends $h_n(y)$ to~$\alpha'_{D,n}$.  From [Mac, (4.2)] we have
$$\prod_{i,j} {1\over 1-x_i y_j} = \sum_\lambda h_\lambda(y) m_\lambda(x).$$
Applying $\theta_{D,y}$ gives
$$\theta_{D,y}\biggl(\prod_{i,j} {1\over 1-x_i y_j}\biggr) =
  \sum_\lambda \alpha_{D,\lambda_1}\alpha_{D,\lambda_2}\cdots m_\lambda(x).$$
Now $\Xi_D(x,0)$ is just the coefficient
of~$a_1a_2\cdots a_d$ in this expression, since this coefficient counts
all path covers of type~$\pi$ exactly $r_\pi!$ times, and
$r_\pi!\, m_\pi(x) = \tilde m_\pi(x)$.  Similarly, $\Xi_{D'}(x,0)$ is the
coefficient of $a_1a_2\cdots a_d$ in
$$\theta'_{D,y}\biggl(\prod_{i,j} {1\over 1-x_i y_j}\biggr) =
  \sum_\lambda \alpha'_{D,\lambda_1}\alpha'_{D,\lambda_2}\cdots m_\lambda(x).$$
Thus it suffices to prove that
$$\omega_x \theta_{D,y}\biggl(\prod_{i,j} {1\over 1-x_i y_j}\biggr)
  = \theta'_{D,y}\biggl(\prod_{i,j} {1\over 1-x_i y_j}\biggr).$$
Now from [Mac, (4.3)] we have
$$ \biggl(\prod_{i,j} {1\over 1-x_i y_j}\biggr) =
   \sum_\lambda s_\lambda(x) s_\lambda(y),$$
so
$$\omega_x \biggl(\prod_{i,j} {1\over 1-x_i y_j}\biggr)
 = \sum_\lambda s_{\lambda'}(x) s_\lambda(y)
 = \sum_\lambda s_\lambda(x) s_{\lambda'}(y)
 = \omega_y \biggl(\prod_{i,j} {1\over 1-x_i y_j}\biggr).$$
Thus 
$$\omega_x \theta_{D,y}
   \biggl(\prod_{i,j} {1\over 1-x_i y_j}\biggr)
  = \theta_{D,y} \omega_x
   \biggl(\prod_{i,j} {1\over 1-x_i y_j}\biggr)
  = \theta_{D,y} \omega_y
   \biggl(\prod_{i,j} {1\over 1-x_i y_j}\biggr).$$
So it suffices to show that $\theta_{D,y}\omega_y = \theta'_{D,y}$.  From
[Mac, (2.6)] we have
$$\sum_n (-1)^n e_n(y) = \biggl(\sum_n h_n(y)\biggr)^{-1},$$
so applying $\theta_{D,y}$ and using \carlitz\ yields
$$\sum_n (-1)^n \theta_{D,y}\bigl(e_n(y)\bigr)
  = \biggl(\sum_n \alpha_{D,n}\biggr)^{-1}
  = \sum_n (-1)^n \alpha'_{D,n}.$$
Equating terms of the same degree, we see that
$$\theta_{D,y}\omega_y\bigl(h_n(y)\bigr) = \theta_{D,y}\bigl(e_n(y)\bigr)
  = \alpha'_{D,n} = \theta'_{D,y}\bigl(h_n(y)\bigr),$$
completing the proof.\qed

Our second proof of \onevarreciprocity\ is similar to
Stanley's proof of the reciprocity theorem for~$\X_G$.
Following Gessel~[Ge2] and Stanley~[St2, section~3], we define a power series
in the variables $x=\{x_1,x_2,\ldots\}$ to be \itc{quasi-symmetric} if
the coefficients of
$$x_{i_1}^{r_1}x_{i_2}^{r_2}\cdots x_{i_k}^{r_k}
  \qquad{\rm and}\qquad
  x_{j_1}^{r_1}x_{j_2}^{r_2}\cdots x_{j_k}^{r_k}$$
are equal whenever $i_1<i_2<\cdots<i_k$ and $j_1<j_2<\cdots<j_k$.
For any subset~$S$ of~$[d-1]$ define the \itc{fundamental}
quasi-symmetric function $Q_{S,d}(x)$ by
$$Q_{S,d}(x) = \sum_{\scriptstyle i_1\le\cdots\le i_d \atop
  \scriptstyle i_j < i_{j+1}\;{\rm if}\; j\in S}
  x_{i_1} x_{i_2} \cdots x_{i_d}.$$
If there is no danger of confusion, we will sometimes
write $Q_S$ for~$Q_{S,d}$ for brevity.

We have the following expansion of $\Xi_D(x;0)$ in terms of fundamental
quasi-symmetric functions.

\proposition\quasisymmetric{If $D$ is a digraph with vertex set~$[d]$,
then
$$\Xi_D(x;0) = \sum_{\pi\in S_d} Q_{S(\pi),d}(x),$$
where $S_d$ is the group of permutations of~$[d]$ and
$$S(\pi) = \{i\in [d] \mid \hbox{$(\pi_i,\pi_{i+1})$ is not an edge
of~$D$} \}.$$}

\proof
We use the expression for $\Xi_D$ given in \coloring.
Given a path coloring of~$D$, arrange the paths in increasing order of
their colors, and within each path arrange the vertices in the order
given by the directed path.  This gives a permutation of the vertices
of~$D$, and it is easy to see that $Q_{S(\pi)}(x)$ counts precisely
the path colorings that give rise to~$\pi$.\qed

We can now give our second proof of \onevarreciprocity.

\noindent{\it Second proof of \onevarreciprocity.}
Without loss of generality we may assume that the vertex set of~$D$
is~$[d]$.  From the same argument as in \quasisymmetric, we see that
$$\Xi_{D'}(x;0) = \sum_{\pi\in S_d} Q_{[d]\backslash S(\pi)}(x).$$
In view of \quasisymmetric, it suffices to show that the map that sends
$Q_S$ to $Q_{[d]\backslash S}$ equals~$\omega$ when restricted to
symmetric functions.  A proof of this fact may be found in the proof
of [St2, Theorem~4.2].\qed

Stanley~[St2] has obtained an analogue of \quasisymmetric\ by
using the theory of acyclic orientations and $P$-partitions.
It is natural to ask if these ideas can be applied to studying~$\Xi_D$.
Unfortunately this does not seem possible.  For example, a key step in
the proof of the analogue of \quasisymmetric\ involves expressing $\X_G$
as a sum of certain poset generating functions, but in general
$\Xi_D$ has no such expression, even if $D$ is acyclic.

\sec Rook theory

As we explained in the previous chapter,
there is a close connection between the path-cycle symmetric function
and rook theory, because $\Xi_D$ is a generalization of the factorial
polynomial of Goldman, Joichi and White.
Thus every theorem about $\Xi_D$ can be
specialized to a theorem in rook theory, and we can also try to
generalize every theorem in rook theory to a theorem about~$\Xi_D$.

A good example of this relationship is \reciprocity, which
can be viewed as a generalization of a result in Riordan~[Rio, Chapter~7,
Theorem~2] relating the
rook numbers of complementary boards, a result which we now state.
If $B$ is a board, we let $B' = \bigl([d]\times[d]\bigr) \backslash B$
denote the complementary board.

\proposition\rookreciprocity{Let $B\subset [d]\times[d]$ be a board.  Then
$R(B';i) = (-1)^d R(B;{}-i-1)$.}

\proof
Let $D$ be the associated digraph.  From
\coverreciprocity\ and \coverrook\ we have
$$R(B';i) = C(D';i,1) = (-1)^d C(D;{}-i-1,1) = (-1)^d R(B;{}-i-1).\qed$$

Riordan's original result is
$$\sum_k r_k^B (d-k)!\,i^k=\sum_k (-1)^k r_k^{B'}\!(d-k)!\,i^k(i+1)^{d-k},$$
which can be shown to be equivalent to \rookreciprocity.  However, it
seems that the formulation of \rookreciprocity, which is much simpler,
has not appeared in the literature before
(although as Gessel~[Ge1] has observed, it follows
immediately from the next proposition below).

Goldman~[Gol] has remarked that \rookreciprocity\ ought to have a direct
combinatorial proof, and indeed it does.  Begin by observing that
$$\eqalign{(-1)^d R(B';{}-i-1) &=
  (-1)^d \sum_{k=0}^d r_k^{B'} ({}-i-1)^{\underline{d-k}} \cr
 &= \sum_{k=0}^d (-1)^k r_k^{B'} (i+d-k)^{\underline{d-k}}.\cr}$$
Now add $i$ extra rows to $[d]\times [d]$.  Then
$r_k^{B'} (i+d-k)^{\underline{d-k}}$ is the number of ways of first
placing $k$~rooks on~$B'$ and then placing $d-k$ more rooks anywhere
(i.e., on $B$, $B'$ or on the extra rows) such that no two rooks
can take each other in the final configuration.
By a straightforward inclusion-exclusion argument,
we see that the resulting configurations in which the set~$S$ of rooks
$B'$ is nonempty cancel out of the above sum, because they are counted
once for each subset of~$S$, with alternating signs.  Thus what survives
is the set of placements of $d$ nontaking rooks on the extended board
such that no rook lies on~$B'$---but it is easy to see that is this
precisely what $R(B;i)$ counts.

We should remark that the idea of adding $i$ extra rows was taken
from~[GJW1].  Also, Gessel~[Ge1] has extended \rookreciprocity\ to the
more general setting of simplicial complexes, and the above combinatorial
proof can be adapted without much difficulty to proving Gessel's more
general result.  We omit the details since they are tangential to our
main purpose.  Finally, we remark that it ought to be possible to
extend the above combinatorial proof to the cover polynomial, but
we have not done so.

We now turn to what is perhaps the most fundamental theorem of rook
theory---the inclusion-exclusion formula (which we have already
mentioned in the proof of \ppositivity)---and generalize it to the
context of the path-cycle symmetric function.
The inclusion-exclusion formula has many equivalent formulations;
the one we shall find most convenient is following one,
whose proof is given implicitly in~[CG1].

\proposition\rookincexcbasic{Let $D$ be a
digraph with $d$~vertices,
and let $N_k^D$ denote the number of ways of placing
$d$~non-taking rooks on $[d]\times[d]$ such that exactly $k$~rooks
lie on the board associated with~$D$.  Then
$$R(D; i) = \sum_k N_k^D {i+k\choose d}.\qed$$}

To state our generalization we need a few more definitions.

\definition
For any pair of integer partitions~$\lambda$ and~$\mu$, define
$D_{\lambda,\mu}$ to be a disjoint union of directed paths and
directed cycles such that the $i$th directed path has $\lambda_i$
vertices and the $j$th directed cycle has $\mu_j$ vertices.
Define
$$\tilde \Xi_{\lambda,\mu} \defeq
   \sum_S {\tilde m_{\pi(S)}(x) \, p_{\sigma(S)}(y) \over \ell(\pi(S))!},$$
where the sum is over all path-cycle covers~$S$ of~$D_{\lambda,\mu}$.
For brevity we shall write $D_\lambda$ for $D_{\lambda,\varemptyset}$
and $\tilde \Xi_\lambda$ for $\tilde \Xi_{\lambda,\varemptyset}$.
We then have the following result.

\theorem\rookincexc{Let $D$ be a digraph with $d$~vertices and let $B$ be the
associated board.  Let ${\Scr N}_{\lambda,\mu}^D$ be the set of placements of
$d$~non-taking
rooks on $[d]\times[d]$ such that the type~$(\pi,\sigma)$ of the path-cycle
cover formed by the set of edges corresponding to rooks placed on~$B$
satisfies ${\rm type}(\pi)=\lambda$ and ${\rm type}(\sigma)=\mu$, and let
$N_{\lambda,\mu}^D = |{\Scr N}_{\lambda,\mu}^D|$.  Then
$$\Xi_D = \sum_{\lambda,\mu} N_{\lambda,\mu}^D \tilde \Xi_{\lambda,\mu},$$
where the sum is over all integer partitions~$\lambda$ and~$\mu$.}

\proof
The proof is similar to the proof of [St3, Theorem~2.3.1].
Given any two integer partitions $\nu$ and~$\eta$,
let ${\Scr R}_{\nu,\eta}^D$ be the set of path-cycle covers~$S$ of~$D$
satisfying ${\rm type}\bigl(\pi(S)\bigr) = \nu$ and
${\rm type}\bigl(\sigma(S)\bigr)=\eta$.
Note that every element of~${\Scr R}_{\nu,\eta}^D$
has $\ell(\nu)$ directed paths
plus some cycles and therefore has a total of $d-\ell(\nu)$ edges.

Now fix a pair of integer partitions $\nu$ and~$\eta$ and
consider the set of pairs $(S,T)$ such that
$S\in {\Scr R}_{\nu,\eta}^D$ and $T$ is an extension of~$S$
(regarded as a placement of non-taking rooks on~$B$)
to a placement of $d$ non-taking rooks on~$[d]\times[d]$.
This set has $\ell(\nu)!\,|{\Scr R}_{\nu,\eta}^D|$ elements,
since for each~$S\in {\Scr R}_{\nu,\eta}^D$
the $\ell(\nu)$ rows and and columns unoccupied by~$S$
can support $\ell(\nu)!$ placements of non-taking rooks.
On the other hand, we can enumerate the set in another way,
by taking each placement of $d$ non-taking rooks on $[d]\times[d]$
and counting how many $S\in{\Scr R}_{\nu,\eta}^D$ it extends.
Now if $T\in {\Scr N}_{\lambda,\mu}^D$ for some $\lambda$ and~$\mu$, then
the number of elements of ${\Scr R}_{\nu,\eta}^D$ that it extends
is just the number~$n_{\lambda,\mu,\nu,\eta}$
of path-cycle covers~$S$ of~$D_{\lambda,\mu}$ satisfying
${\rm type}\bigl(\pi(S)\bigr) = \nu$ and
${\rm type}\bigl(\sigma(S)\bigr)=\eta$.  Since every placement of
$d$~non-taking rooks on $[d]\times[d]$ belongs to ${\Scr N}_{\lambda,\mu}^D$
for some $\lambda$ and~$\mu$, we have
$$\sum_{\lambda,\mu} N_{\lambda,\mu}^D n_{\lambda,\mu,\nu,\eta} =
   \ell(\nu)!\,|{\Scr R}_{\nu,\eta}^D|.$$
Now divide both sides by $\ell(\nu)!$,
multiply both sides by $\tilde m_\nu(x)\,p_\eta(y)$, and sum
over all $\nu$ and~$\eta$ to obtain the desired result.\qed

The next proposition shows that \rookincexc\ does indeed generalize
\rookincexcbasic, a fact which may not be obvious at first glance.

\proposition\rookspecialize{For any integer partitions $\lambda$ and~$\mu$,
$$\tilde \Xi_{\lambda,\mu}(1^i;1) = {i+d-\ell(\lambda) \choose d},$$
where $d = |\lambda|+|\mu|$.}

\proof
Directly from the definitions we have
$$\tilde \Xi_{\lambda,\mu}(1^i;1) =
  \sum_S {i^{\underline{\ell(\pi(S))}} \over \ell(\pi(S))!}
 =\sum_S {i \choose \ell(\pi(S))},$$
where the sum is over all path-cycle covers of~$D_{\lambda,\mu}$.
The sum can be broken up into a double sum:
$$\tilde \Xi_{\lambda,\mu}(1^i;1) =
  \sum_k \sum_{\{S\mid \ell(\pi(S))=k\}} {i \choose k}.$$
But $\ell(\pi(S)) = k$ if and only if $|S| = d - k$.
Since every subset of the edges of~$D_{\lambda,\mu}$ is a path-cycle
cover, and since the total
number of edges of~$D_{\lambda,\mu}$ is $d-\ell(\lambda)$, we have
$$\tilde \Xi_{\lambda,\mu}(1^i;1) =
  \sum_k {d-\ell(\lambda) \choose d - k} {i \choose k} = 
  {i + d - \ell(\lambda) \choose d}.\qed$$

In view of \coverrook, \rookspecialize, and the fact that
the number of edges of a path-cycle cover of type~$(\pi,\sigma)$ is
$d-\ell(\pi)$, we see that \rookincexc\ implies \rookincexcbasic.

We remark that in passing that Gessel~[Ge1] has obtained a generalization
of \rookincexcbasic\ for the cover polynomial that does not appear
to follow from our results.

It might seem that \rookincexc\ is contrived, since we seem to have
defined $\tilde \Xi_{\lambda,\mu}$ just so that \rookincexc\ would come out
right.  In fact, however, the functions $\tilde \Xi_{\lambda,\mu}$
are surprisingly interesting objects in their own right.
For a start, we have the following easy fact.

\proposition\xibasis{The functions $\tilde \Xi_\lambda$
form a linear basis for the ring of symmetric functions over the rationals,
and the functions $\tilde \Xi_{\lambda,\mu}$ form a linear basis
for the ring of symmetric functions in two sets of variables (again over the
rationals).}

\proof
Write $$\tilde \Xi_\lambda = \sum_\mu c_{\lambda,\mu} \tilde m_\mu.$$
Then it is clear from the definition of~$\tilde \Xi_\lambda$ that
$c_{\lambda,\lambda}\ne0$ and $c_{\lambda,\mu}\ne0$ only if
$\lambda\ge\mu$ in refinement order.  Thus the matrix $(c_{\lambda,\mu})$
with respect to any linear extension of refinement order is triangular
with nonzero entries on the diagonal.  This proves the first assertion.
To prove the second assertion, define a partial order on pairs of
integer partitions by setting $(\lambda,\mu)<(\nu,\eta)$ if the multiset
of parts of~$\eta$ can be partitioned into two multisets $\alpha$
and~$\beta$ such that $\mu=\beta$ and $\lambda$ is a refinement of
$\nu\cup\alpha$.  The same kind of reasoning as before, with this partial
order in place of refinement order and with the basis
$\tilde m_\lambda(x)\,p_\mu(y)$ in place of $\tilde m_\lambda$,
can then be applied to prove the second assertion.\qed

The $\tilde\Xi_\lambda$ turn out to be particularly interesting, as we
shall see presently.  The following table expresses $\tilde\Xi_\lambda$
in terms of the monomial symmetric functions for some small values
of~$\lambda$.
{
\def\normalbaselines{\baselineskip=14pt}
\def\vp{\vphantom{2\over3}}
$$
\eqalignno{
\left[\matrix{\tilde\Xi_2\hfill\cr \tilde\Xi_{11}\hfill\cr}\right] &=
\left[\matrix{1&1\cr 0&1\cr}\right]
\left[\matrix{m_2\hfill\cr m_{11}\hfill\cr}\right]\cr
\left[\matrix{\tilde\Xi_3\hfill\cr \tilde\Xi_{21}\hfill\cr
 \tilde\Xi_{111}\hfill\cr}\right] &=
\left[\matrix{1&1&1\cr 0&{1\over2}&1\cr 0&0&1\cr}\right]
\left[\matrix{m_3\hfill\cr m_{21}\hfill\cr m_{111}\hfill\cr}\right]\cr
\left[\matrix{\tilde\Xi_4\hfill\cr \tilde\Xi_{31}\hfill\cr
 \tilde\Xi_{22}\hfill\cr \tilde\Xi_{211}\hfill\cr
 \tilde\Xi_{1111}\hfill\cr} \right]&=
\left[\matrix{1&1&1&1&1\cr 0&{1\over 2}&0&{2\over 3}&1\cr
 0&0&1&{2\over3}&1\cr 0&0&0&{1\over3}&1\cr 0&0&0&0&1\cr}\right]
\left[\matrix{m_4\hfill\cr m_{31}\hfill\cr m_{22}\hfill\cr m_{211}\hfill
 \cr m_{1111}\hfill\cr}\right]\cr
\left[\matrix{\tilde\Xi_5\hfill\cr \tilde\Xi_{41}\hfill\cr
 \tilde\Xi_{32}\vp\hfill\cr \tilde\Xi_{311}\vp\hfill\cr
 \tilde\Xi_{221}\hfill\vp\cr \tilde\Xi_{2111}\hfill\cr
 \tilde\Xi_{11111}\hfill\cr}\right] &=
\left[\matrix{1&1&1&1&1&1&1\cr 0&{1\over2}&0&{2\over3}&{1\over3}&{3\over4}&1\cr
 0&0&{1\over2}&{1\over3}&{2\over3}&{3\over4}&1\cr
 0&0&0&{1\over3}&0&{1\over2}&1\cr 0&0&0&0&{1\over3}&{1\over2}&1\cr
 0&0&0&0&0&{1\over4}&1\cr 0&0&0&0&0&0&1\cr}\right]
\left[\matrix{m_5\vp\hfill\cr m_{41}\vp\hfill\cr m_{32}\vp\hfill\cr
 m_{311}\hfill\vp\cr m_{221}\vp\hfill\cr m_{2111}\vp\hfill\cr
 m_{11111}\vp\hfill\cr}\right]\cr
}$$
}

While a casual inspection of the above table may not reveal any
interesting patterns, we have the following surprising result.

\theorem\involution{The linear map that sends $\tilde\Xi_\lambda$
to $(\sgn\lambda)\tilde m_\lambda/\ell(\lambda)!$ is an involution.}

\proof
Given any integer partitions $\lambda$ and~$\mu$, let $\pi$ be any set
partition of type~$\lambda$ and define
$$c_{\lambda,\mu} \defeq
  \sum_{\{\sigma\ge\pi\mid {\rm type}(\sigma)=\mu\}}
  \lambda(\pi,\sigma)!.$$
Note that $c_{\lambda,\mu}$ does not depend on the choice of~$\pi$.
We claim that
$$\Xi_{D_\mu} =
  \sum_\lambda {r_\mu!\over r_\lambda!} c_{\lambda,\mu} \tilde m_\lambda.$$
To see this, first consider the case where $r_\mu!=r_\lambda!=1$, i.e.,
the case of distinct parts.  We have a disjoint union~$D_\mu$
of directed paths and
we want to count the number of path covers of type~$\lambda$.  In a path
cover of~$D_\mu$, each directed path is broken up into a sequence of
smaller directed paths.  So the path covers can be enumerated as follows:
take a set partition~$\pi$ of type~$\lambda$ and consider all
ways of grouping its blocks into a partition~$\sigma$ of type~$\mu$
and then linearly ordering the blocks of~$\pi$ within each block
of~$\sigma$.  Such a configuration determines a path cover: for any
block~$b$ of~$\sigma$, the sequence of blocks of~$\pi$ in~$b$
dictates the sizes of the sequence of smaller directed paths
composing the directed path in~$D_\mu$ corresponding to~$b$.
It is easy to see that this correspondence is bijective, and this
proves our claim in the case of distinct parts.  For the general
case, observe that we want equal-sized parts of~$\pi$ to be
indistinguishable and equal-sized parts of~$\sigma$ to be distinguishable,
so we must multiply by $r_\mu!/r_\lambda!$.

\relax From \doubilet\ we see that the matrix
$\bigl((\sgn\lambda)c_{\lambda,\mu}\bigr)$ is the matrix of~$\omega$
(with respect to the augmented monomial symmetric function basis)
and is therefore an involution.  From our claim it follows that the
matrix relating
$(\sgn\lambda)\tilde m_\lambda / r_\lambda!$ and
$\Xi_{D_\mu}/r_\mu!$
or equivalently the matrix relating
$${(\sgn\lambda)\tilde m_\lambda \over r_\lambda!\ell(\lambda)!}
  \qquad{\rm and}\qquad {\tilde \Xi_\mu\over r_\mu!}$$
is an involution.  But then the desired result follows, since the
factors of $r_\lambda!$ and~$r_\mu!$ amount to conjugating by a
(diagonal) matrix, and this does not change the involution property.\qed

Notice the close connection between the involution of \involution\ and
the involution~$\omega$.  (In fact, it was a suggestion by Stanley that
the two involutions might be equal that led to the proof of \involution.)
The two involutions are not the same, however---the former is not even a
homomorphism---and philosophically speaking it is still unclear why they
are related.  For example, if we specialize to the polynomial level,
$\omega$ becomes (essentially) the operation of substituting $-x$ for~$x$,
but we do not know of any such intuitive interpretation for the involution
of~\involution.

The $\tilde \Xi_\lambda$ are also closely related to
the fundamental quasi-symmetric functions~$Q_S$ defined
in the previous section.
(The knowledgeable reader may already have suspected this
since the fundamental quasi-symmetric functions specialize
to the same polynomial basis as the $\tilde\Xi_\lambda$ do.)
If $S$ is a subset of~$[d-1]$, then we define the \itc{type}
of~$S$ to be the integer partition whose parts are the lengths
of the subwords obtained by breaking the word $123\ldots d$
after each element of~$S$.

\theorem\xitquasi{Let $g$ be any symmetric function.  If $a_\lambda$
and $b_S$ are constants such that
$$g = \sum_\lambda a_\lambda \tilde \Xi_\lambda
 \qquad{\rm and}\qquad g = \sum_S b_S Q_S,$$
then
$$ a_\lambda = \sum_{\{S \mid {\rm type}(S)=\lambda\}} b_S.$$}

\proof
It is not difficult to see that
it suffices to prove the theorem for the case $g=\tilde\Xi_\mu$.
Let $d$ be any positive integer and let $S$ be any subset of $[d-1]$
and define
$$\tilde Q_S \defeq \sum_{\scriptstyle i_1 \le i_2 \le \cdots \le i_d
  \atop \scriptstyle i_j < i_{j+1}\ {\rm iff}\ j \in S}
  x_{i_1} x_{i_2} \cdots x_{i_d}.$$
Then
$$m_\lambda = \sum_{\{S\mid {\rm type}(S) = \lambda\}} \tilde Q_S
  \qquad{\rm and}\qquad
  Q_S = \sum_{T\supset S} \tilde Q_T.$$
By an inclusion-exclusion argument,
$$m_\lambda = \sum_{\{S\mid {\rm type}(S) = \lambda\}} \;
   \sum_{T\supset S} (-1)^{|T| - |S|} Q_T.$$
Let $q_{\lambda T}$ be the coefficient of~$Q_T$ in~$m_\lambda$.
We compute
$$ \sum_{\{ T\mid {\rm type}(T) = \nu\}} q_{\lambda T}.$$
Observe that there is a bijection
between subsets of type~$\lambda$
and orderings of the parts of~$\lambda$: given a subset~$S\subset[d-1]$
of type~$\lambda$, take the sequence of the lengths of the subwords
of the word $123\ldots d$ obtained by breaking after each
element of~$S$.
Thinking of such subwords as directed paths, we see that for any
fixed $S$ of type~$\lambda$, the number of subsets $T\supset S$
such that ${\rm type}(T) = \nu$ is just the number of path-coverings
of $D_\lambda$ of type~$\nu$, which from the proof of \involution\ is
$${r_\lambda!\over r_\nu!} \,c_{\nu\lambda}$$
(using the notation of \involution).  Now there are
$\ell(\lambda)!/r_\lambda!$ subsets~$S$ of type~$\lambda$,
and if ${\rm type}(S)=\lambda$ and ${\rm type}(T)=\nu$ then
$$(-1)^{|T|-|S|} = ({\rm sgn}\,\nu)({\rm sgn}\,\lambda).$$
Putting all this together, we see that
$$ \sum_{\{ T\mid {\rm type}(T) = \nu\}} \!\!\!q_{\lambda T}
  = {\ell(\lambda)! \over r_\nu!} \,c_{\nu\lambda} 
    ({\rm sgn}\,\nu)({\rm sgn}\,\lambda).$$
But, again from the proof of \involution,
$$\tilde \Xi_\mu = \sum_\lambda {r_\mu! \over r_\lambda!} \,c_{\lambda\mu}
   {r_\lambda!\over \ell(\lambda)!} m_\lambda.$$
Hence if $g=\tilde \Xi_\mu$, then
$$\eqalign{ \sum_{\{S\mid {\rm type}(S) = \nu\}} \!\!\!b_S &= 
    \sum_\lambda {r_\mu! \over r_\lambda!} \,c_{\lambda\mu}
    {r_\lambda!\over \ell(\lambda)!} \cdot
    {\ell(\lambda)!\over r_\nu!} \,c_{\nu\lambda}
    ({\rm sgn}\,\nu)({\rm sgn}\,\lambda)\cr
 &= {r_\mu!\over r_\nu!} \sum_\lambda 
    ({\rm sgn}\,\nu)c_{\nu\lambda}({\rm sgn}\,\lambda)c_{\lambda\mu}\cr
 &= \delta_{\mu\nu},\cr}$$
because $\bigl(({\rm sgn}\,\lambda)c_{\lambda\mu}\bigr)$ is the matrix
of~$\omega$ with respect to the augmented
monomial symmetric function basis (by \doubilet), and
$\omega$ is an involution.  This completes the proof.\qed

\corollary\qpositive{If $g$ is a $Q$-positive symmetric function,
then $g$ is $\tilde \Xi$-positive.  In particular, $\X_G$ and
the Schur functions are $\tilde \Xi$-positive.}

\proof The first assertion is immediate from \xitquasi.
The fact that $\X_G$ and the Schur functions are $Q$-positive is
``folklore''; it is implicit in~[St5], but see also [Ge2] and~[St2].\qed

We caution the reader not to read more into \xitquasi\ than is actually
there!  For example, $\tilde\Xi_\lambda$ is \itc{not} $Q$-positive.
Nor is it true that the only $Q_T$'s in the
$Q$-expansion of~$\tilde\Xi_\lambda$ with nonzero coefficients
are those with ${\rm type}(T)=\lambda$.
Thus, while \xitquasi\ allows one to translate combinatorial
\itc{interpretations of the coefficients} of the $Q$-expansion
of a symmetric function~$g$ into combinatorial interpretations of
the the coefficients of the $\tilde\Xi$-expansion of~$g$, there is
no guarantee that combinatorial \itc{proofs} can be so translated.
Some tricky reshuffling of combinatorial information occurs in the
transition from the $Q$'s to the~$\tilde\Xi$'s.

We should mention another, somewhat more philosophical, reason that
the $\tilde\Xi_\lambda$ are interesting.
The bases $\tilde m_\lambda$, $p_\lambda$, $e_\lambda$, $h_\lambda$,
$s_\lambda$, and $f_\lambda$ occur frequently ``in nature.''  Similarly,
there are certain ``natural'' bases for polynomials, and moreover there
is a correspondence between some of the symmetric function bases and
the polynomial bases given by $g\mapsto g(1^i)$, e.g.,
$p_\lambda$ corresponds to~$i^n$, and the reciprocally
related bases $\tilde m_\lambda$ and $f_\lambda$ correspond to
$i^{\underline n}$ and~$i^{\overline n}$.
However, so far as we are aware, no symmetric function counterpart
to the polynomial basis
$${i+n\choose d}_{n=0,1,\ldots,d}$$
has been proposed before.  \rookspecialize, together
with \involution\ and \qpositive, suggests that the
$\tilde\Xi_\lambda$ may be the ``right'' symmetric function
generalization of this polynomial basis.

If this is the case, one might hope that specializing
$\tilde\Xi_{\lambda,\mu}$ might give rise to an interesting
basis for polynomials in two variables.
Unfortunately, this does not seem to be true.
However, specializing $\tilde\Xi_{\lambda,\mu}$ to a polynomial does
provide a simple proof of a theorem of Chung and Graham
whose original proof is quite complicated.
Following Chung and Graham,
for any placement~$T$ of $d$ non-taking rooks on $[d]\times[d]$,
let ${\rm drop}(T)$ be the subgraph of~$D$ with edges corresponding
to the squares occupied by the rooks of~$T$.
If $D$ is a digraph with $d$ vertices,
let $\delta_D(q,r,s)$ be the number of ordered pairs $(S,T)$
such that $S$ is a set of $r$~edges of~$D$ forming precisely $s$ disjoint
cycles and $T$ is a placement of $d$ non-taking
rooks on $[d]\times[d]$ with $S\subset {\rm drop}(T)$
and $|{\rm drop}(T)| = q+r$.
Chung and Graham's result ([CG1, Theorem~2] or [CG2, Theorem~3])
is then the following.

\proposition\chunggraham{For any digraph~$D$ with $d$ vertices,
$$C(D; i,j) = \sum_{q,r,s} \delta_D(q,r,s) {i+q \choose d-r} (j-1)^s.$$}

\proof
We can restate the desired result as
$$\eqalign{C(D; i,j+1) &= \sum_{q,r,s,t} \delta_D(q,r,s) {q \choose d-r-t}
   {i \choose t} j^s \cr
 &= \sum_{q,r,s,t}
    \delta_D(q,r,s) {q \choose t-d+r+q} {i\choose t} j^s.\cr}$$
\relax From \rookincexc\ and the definition of~$\tilde \Xi_{\lambda,\mu}$
we have
$$C(D; i,j+1) = \sum_{\lambda,\mu,t,u} N_{\lambda,\mu}^D
  n_{\lambda,\mu,t,u} {i\choose t} (j+1)^u
 = \sum_{\lambda,\mu,s,t,u} N_{\lambda,\mu}^D n_{\lambda,\mu,t,u}
   {u \choose s} {i \choose t} j^s,$$
where $n_{\lambda,\mu,t,u}$ is the number of path-cycle covers
of~$D_{\lambda,\mu}$ with $t$~paths and $u$~cycles.  Now ${i\choose t} j^s$
is a basis for polynomials in two variables, so equating coefficients we see
that we just need to prove that for any fixed $s$ and~$t$,
$$\sum_{\lambda,\mu,u} N_{\lambda,\mu}^D n_{\lambda,\mu,t,u} {u \choose s}
 = \sum_{q,r} \delta_D(q,r,s) {q \choose t-d+r+q}.$$
We can think of both sides as counting placements of $d$ non-taking
rooks on $[d]\times[d]$ with certain multiplicities.  On the left-hand
side, the number of times each such placement~$T$ is counted equals the
number of path-cycle covers of ${\rm drop}(T)$ with exactly $t$~paths
plus some number of cycles of which $s$ are distinguished.
As for the right-hand side, we can rewrite it as
$$\sum_{e,r} \delta_D(e-r,r,s) {e-r \choose t-d+e}.$$
Then for any placement~$T$, only one value of~$e$ (namely,
$e=|{\rm drop}(T)|$) involves~$T$.
Thus if we let $e=|{\rm drop}(T)|$,
the number of times $T$ is counted is
$$\sum_r {\hbox{the number of ways of choosing} \strut\atopwithdelims()
  \hbox{$s$ cycles of ${\rm drop}(T)$ with $r$ edges}}
  \cdot {e-r \choose t-d+e},$$
which is just the number of ways of choosing $s$ cycles and then deleting
$t-(d-e)$ of the remaining edges (i.e., creating $t-(d-e)$ new paths).
But $d-e$ is the original number of paths in ${\rm drop}(T)$, so this
results in a total of exactly $t$~paths.  The proposition follows.\qed

Let us now return from this digression to the problem of generalizing
rook theory to the context of the path-cycle symmetric function.
Our next result generalizes a M\"obius inversion formula for factorial
polynomials due to Goldman, Joichi and White~[GJW4].
For simplicity we consider only the case of acyclic digraphs,
although the generalization to arbitrary digraphs is straightforward.
So suppose $D$ is an acyclic digraph with $d$ vertices
and let $B$ be its associated board.
Following an idea of Goldman, Joichi and White,
extend the columns of~$[d]\times[d]$ infinitely
downwards, so that there are now infinitely many rows.
Let $\Scr S$ be the set of all placements of $d$~rooks such that

\item{1.} every rook lies either on~$B$ or one of the appended squares, and

\item{2.} no two rooks lie in the same column.

Given $S\in{\Scr S}$, define $\pi(S)$ to be the partition of~$[d]$ in which
two numbers $i$ and~$j$ lie in the same block if and only if the rooks in
columns $i$ and~$j$ lie in the same row.  To each $S\in{\Scr S}$ we also
associate a coloring of~$[d]$ as follows.
Color the vertex $i\in[d]$ with color~$j$
if the rook in column~$i$ lies in the $j$th \itc{appended} row.
Otherwise, if the rook in column~$i$ lies in the $j$th \itc{original} row,
make vertex~$i$ the same color as vertex~$j$.  Since there is exactly one
rook in each column, and since $D$ is acyclic,
these rules give a well-defined coloring~$c_S$.
For every set partition of~$[d]$, define
$$T_\pi^D = T_\pi^B \defeq \sum_{\{S\in{\Scr S}\mid \pi(S)=\pi\}} x^S,$$
where
$$x^S \defeq \prod_{i\in[d]} x_{c_S(i)}.$$
Finally define
$$T_{\ge\pi}^D = T_{\ge\pi}^B \defeq \sum_{\sigma\ge\pi} T_\sigma^B.$$

\proposition\gjw{For any acyclic digraph~$D$ with $d$ vertices,
$$\Xi_D = \sum_{\pi\in\Pi_d} \mu(\hat 0,\pi) T_{\ge\pi}^D.$$}

\proof
By M\"obius inversion, the right-hand side is just~$T_{\hat 0}^D$.
The $S\in{\Scr S}$ such that $\pi(S)=\hat 0$ are just the placements
in which no two rooks lie in the same row or column.  The rooks on~$B$
then define a path cover and the rooks on the appended rows then
ensure that distinct paths are assigned distinct colors.
The theorem follows from \coloring.\qed

It is not hard to show that this result specializes to~[GJW4, Theorem~1(a)].
One might
again object that \gjw\ is contrived because $T_{\ge\pi}^D$ is simply
a formal device to represent what one gets by M\"obius inversion.  This
time the objection is harder to meet, because $T_{\ge\pi}^D$ is not as
``nice'' an object as~$\tilde \Xi_{\lambda,\mu}$.  For example,
${\rm type}(\pi)={\rm type}(\sigma)$ does \itc{not} imply
$T_{\ge\pi}^D = T_{\ge\sigma}^D$.  However, we do have
one result that gives some more information about~$T_{\ge\pi}^D$.

\proposition\ppositive{If $D$ is an acyclic digraph then
$T_{\ge\pi}^D$ is $p$-positive.}

\proof
Let $B$ be the associated board.  We have
$$T_{\ge\pi}^D = \sum_{\{S\in{\Scr S}\mid \pi(S)\ge\pi\}} x^S.$$
Collect terms that have identical placements of rooks on~$B$.  From
the definitions we see that each such collection of terms
corresponds to the set of colorings of~$V(D)$ that
are monochromatic on the connected components of the subgraph
of~$D$ whose edges are those selected by the placement of rooks on~$B$,
except that the condition $\pi(S)\ge\pi$ imposes the further condition
that components which contain elements of the same block of~$\pi$ must
always be colored the same color.  This gives a power sum symmetric
function, so $T_{\ge\pi}^D$ is a sum of power sums.\qed

Note that while \gjw\ resembles Stanley's formula [St2, Theorem~2.6]
$$\X_G = \sum_{\pi\in L_G} \mu(\hat 0,\pi) p_\pi$$
(where $L_G$ is the lattice of contractions of~$G$), there is a significant
difference in that in \gjw\ all the dependence on the digraph is contained
in the~$T_{\ge\pi}^D$ whereas for~$\X_G$ all the dependence is
contained in~$L_G$.  A variant of of \gjw\ can be obtained by
considering rook placements with no two rooks in the same
\itc{row}, but this result also seems contrived and does not
suggest any satisfactory analogue of~$L_G$, so we omit the details.

More rook theory can undoubtedly be generalized to our symmetric
function context, but we shall now turn to a different aspect of~$\Xi_D$.

\sec The Poset Chain Conjecture

One of the original
motivations for studying $\X_G$ and~$\Xi_D$ is a conjecture by
Stanley and Stembridge~[S-S, Conjecture~5.5]
called the \itc{Poset Chain Conjecture.}
We restate this conjecture here for convenience.
Following Stanley~[St2, section~5],
we write $\bf a+b$ for the poset that is a disjoint
union of an $a$-element chain and a $b$-element chain, and we say that a
poset is \itc{$\bf (a+b)$-free} if it contains no induced subposet isomorphic
to $\bf a+b$.
Then the Stanley-Stembridge conjecture is equivalent to the following.

\conjecture\epositivity{If $P$ is a $\bf (3+1)$-free poset, then
$\X_{G(P)}$ is $e$-positive.}

In view of \connection, this conjecture can also be viewed as a conjecture
about~$\Xi_D$.  One of the most important partial results is the following
theorem of Gasharov~[Ga1].

\proposition\gasharov{If $P$ is a $\bf (3+1)$-free poset, then
$\X_{G(P)}$ is $s$-positive.\qed}

The main result of this section is
a slight extension of Gasharov's result that will illustrate
the subtlety of \epositivity.  To state our result we need some definitions.

\definition
A loopless digraph is \itc{weakly $\bf (3+1)$-free} if, for any ordered
pair $(u,v)$ of vertices of~$D$, either $D$ or~$D'$ fails to have a directed
path of length two from $u$ to~$v$.

Note that weakly $\bf (3+1)$-free digraphs need not be transitively closed
or even acyclic.  Our nomenclature is justified by the following proposition.

\proposition\weaklyfree{If $P$ is a poset, then $P$ is $(\bf 3+1)$-free
if and only if $D(P)$ is weakly $\bf (3+1)$-free.}

\proof
Saying that $D(P)$ is weakly $\bf (3+1)$-free is equivalent to saying that
if $u\to v\to w$ is a directed path of length two in~$D(P)$
and $x$ is any element such that $(u,x)$ is \itc{not} an edge of~$D(P)$,
then $(x,w)$ is an edge of~$D(P)$.
Saying that $P$ is $(\bf 3+1)$-free is equivalent to saying that
if $u\to v\to w$ is a chain in~$P$ and $x$ is any element such that
$u\not<x$ in~$P$, then $x<w$ in~$P$.  Clearly these two are equivalent.\qed

\definition
Let $D$ be a digraph.  A \itc{$D$-array} is an array
$$\matrix{v_{1,1}&v_{1,2}&\ldots\cr
          v_{2,1}&v_{2,2}&\ldots\cr
          \ldots&&\cr}$$
where each $v_{i,j}$ is either undefined or an element of~$D$ and such
that

\item{1.} for all $i,j\ge1$, if $v_{i,j+1}$ is defined,
then $v_{i,j}$ is defined and
$(v_{i,j}, v_{i,j+1})$ is an edge of~$D$, and

\item{2.} every element of~$D$ appears exactly once in the array.

The \itc{shape} of a $D$-array is the sequence of
the lengths of (the defined portion of) the rows.
A \itc{$D$-tableau} is a $D$-array such that

\item{3.} for all $i,j\ge1$, if $v_{i+1,j}$ is defined,
then $v_{i,j}$ is defined and
$(v_{i+1,j},v_{i,j})$ is \itc{not} an edge of~$D$.

Our definitions of $D$-array and $D$-tableau are motivated by Gasharov's use
of Gessel-Viennot~[GV2] $P$-arrays and $P$-tableaux in his proof of \gasharov.
We can now state our generalization.

\theorem\spositivity{If $D$ is a weakly $\bf (3+1)$-free digraph, then
the coefficient of~$s_\lambda$ in $\Xi_D(x,0)$ is the number of $D$-tableaux
of shape~$\lambda$.}

\proof
The proof is almost identical to Gasharov's,
and we refer to his paper for some details which we shall omit.
Let $S_\ell$ denote the group of permutations of~$[\ell]$.
If $\lambda=(\lambda_1,\ldots,\lambda_\ell)$ is an integer partition
and $\pi\in S_\ell$, then we denote by $\pi(\lambda)$ the sequence
$$(\lambda_{\pi(j)} - \pi(j) + j)_{j=1}^\ell.$$
Define $c_\lambda$ by
$$\Xi_D(x,0) = \sum_\lambda c_\lambda s_\lambda(x).$$
By the same Jacobi-Trudi argument that Gasharov uses,
$$c_\lambda = \sum_{\pi\in S_\ell} (\sgn\pi) \cdot
  \bigl(\hbox{coefficient of $\prod x_i^{\pi(\lambda)_i}$
     in $\Xi_D(x,0)$}\bigr),$$
where $\sgn\pi$ is the sign of the permutation~$\pi$.
Now by \coloring, $\Xi_D(x,0)$ counts path colorings of~$D$, and path
colorings of~$D$ are in bijection with $D$-arrays (the rows of the $D$-array
give the directed paths and the path in row~$i$ is assigned the color~$i$).
If we let
$$A = \{(\pi,T) \mid \hbox{$\pi\in S_\ell$ and $T$ is a $D$-array of
           shape $\pi(\lambda)$} \},$$
it then follows that
$$c_\lambda = \sum_{(\pi,T) \in A} \sgn\pi.$$
Now let
$$B = \{(\pi,T)\in A \mid \hbox{$T$ is \itc{not} a $D$-tableau}\}$$
and note that if $T$ is a $D$-tableau, then $\pi(\lambda)_1\ge
\pi(\lambda)_2\ge\cdots$ so that $\pi$ must be the identity permutation.
Thus to prove the theorem it suffices to find an involution
$\varphi:B \to B$ such that if $(\sigma,T')=\varphi(\pi,T)$ then
$\sgn\sigma=-\sgn\pi$.  Gasharov's involution works without modification;
for completeness we restate it here.  If
$$T = \matrix{v_{1,1}&v_{1,2}&\ldots\cr
          v_{2,1}&v_{2,2}&\ldots\cr
          \ldots&&\cr}$$
then let $c=c(T)$ be the smallest positive integer such that condition~3
fails for $j=c$ and some~$i$.  Let $r=r(T)$ be the largest $i$ with this
property.  Define $\sigma = \pi \circ (r,r+1)$ where $(r,r+1)$ is the
permutation that interchanges $r$ and~$r+1$.  Define
$$T' = \matrix{u_{1,1}&u_{1,2}&\ldots\cr
          u_{2,1}&u_{2,2}&\ldots\cr
          \ldots&&\cr}$$
by letting

\item{(a)} $u_{i,j} = v_{i,j}$ if $i\ne r$ or $i\ne r+1$ or ($i=r$ and
$j\le c-1$) or ($i=r+1$ and $j\le c$);
\item{(b)} $u_{r,j} = v_{r+1,j+1}$ if $j\ge c$ and $v_{r+1,j+1}$ is defined;
\item{(c)} $u_{r+1,j} = v_{r,j-1}$ if $j\ge c+1$ and $v_{r,j-1}$ is defined.

(Other values of the array~$T'$ remain undefined.)
Now row $r+1$ of~$T'$ satisfies condition~1, because
if $v_{r,c}$ is defined then $(v_{r+1,c},v_{r,c})$ is an edge of~$D$ by
definition of $r$ and~$c$.
To show that $T'$ is a $D$-array it suffices 
to show that row~$r$ satisfies condition~1 since condition~2 is obviously
satisfied.  Possible trouble arises only if $c\ge2$, but then
$v_{r+1,c-1} \to v_{r+1,c} \to v_{r+1,c+1}$ is a path
of length two in~$D$ and $(v_{r+1,c-1}, v_{r,c-1})$ is \itc{not} an edge
in~$D$, so it follows from the assumption that $D$ is weakly $(\bf 3+1)$-free
that $(v_{r,c-1},v_{r+1,c+1})$ is an edge of~$D$, and condition~1 is met.
Now $(v_{r+1,c},v_{r+1,c+1})$ is an edge of~$D$ so if $u_{r,c}$ is defined
$(u_{r+1,c},u_{r,c})$ is an edge of~$D$ (since $u_{r+1,c}=v_{r+1,c}$ and
$u_{r,c}=v_{r+1,c+1}$) and thus $T'$ is not a $D$-tableau.  It is clear
that $T'$ has shape~$\sigma(\lambda)$ and that $c(T')=c(T)$ and
$r(T')=r(T)$.  Also, $\sgn\sigma=-\sgn\pi$,
so $\varphi$ is the desired sign-reversing involution.\qed

Note that a digraph is weakly $(\bf 3+1)$-free if and only
if its complement is weakly $(\bf 3+1)$-free, so \onevarreciprocity\
applied to \spositivity\ does not enlarge the class of known $s$-positive
path-cycle symmetric functions.

It is natural to conjecture that if $D$ is weakly $(\bf 3+1)$-free then
$\Xi_D(x,0)$ is $e$-positive, but for instance if we let $D$ be
the  digraph

\vskip 0.3truein
\hskip 2truein
\psfig{figure=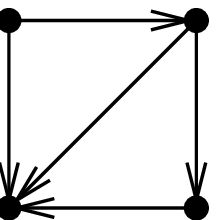}

\noindent
we find (with the aid of John Stembridge's SF package for Maple) that
$$\Xi_D(x,0) = s_4 + 2s_{31} + s_{22} + 4s_{211} + 3s_{1111}
   = 3e_{31} - e_{211} + e_{1111}.$$
In fact, of the five essentially distinct weakly $(\bf 3+1)$-free acyclic
digraphs on four vertices that are not transitively closed, only one is
$e$-positive.  So the way the property of being $(\bf 3+1)$-free is
used in Gasharov's proof is far from enough to yield $e$-positivity
even if the condition of acyclicity is added.  This shows how delicate
\epositivity~is.


\chap The Chromatic Symmetric Function

In this chapter we prove a number of miscellaneous results about the
chromatic symmetric function~$\X_G$.  The basic reference for~$\X_G$
is [St2] (but see also [St4][Ga1][Ga2]).

\sec $G$-ascents

As explained in~[St2],
the expansion of~$\X_G$ in terms of fundamental quasi-symmetric
functions has an interpretation in terms of $P$-partitions.
In the previous chapter we saw that this implies
that the coefficients of the $\tilde\Xi$-expansion of~$\X_G$
have a combinatorial interpretation in terms of $P$-partitions.
In this section we give another
combinatorial interpretation of the coefficients that is based
on the concept of a $G$-descent (see~[CG2]; the definition in~[CG1]
contains a minor error) or the equivalent concept of a $G$-ascent.

\definition
Fix a graph~$G$ with vertex set~$[n]$.  Given a permutation $\pi$ of~$[n]$
and a vertex $v\in[n]$,
define the \itc{rank~$\rho_\pi(v)$} of~$v$ to be the largest integer~$r$
for which there exists an increasing sequence of positive integers
$$i_1 < i_2 < \cdots < i_r < i_{r+1} = \pi^{-1}(v)$$
such that $\{\pi(i_j), \pi(i_{j+1})\}$ is an edge of~$G$ for all~$j$.
We say that $\pi$ has a \itc{$G$-ascent} at~$v$ if either

\item{(i)} $\rho_\pi(v) < \rho_\pi(w)$, or

\item{(ii)} $\rho_\pi(v) = \rho_\pi(w)$ and $v<w$,

\noindent where $w = \pi\bigl(\pi^{-1}(v) + 1\bigr)$.
The \itc{$G$-ascent type} of~$\pi$ is the integer partition of~$n$
whose parts are the lengths of the subwords obtained by breaking the
one-line representation of~$\pi$ after each number at which $\pi$ has
a $G$-ascent.

For example, suppose that $n=3$ and that $G$ has an edge between $1$ and~$3$
and no other edges.  Let $\pi$ be the permutation~$132$, i.e., the
permutation that fixes~1 and exchanges 2 and~3.  Then
$$\rho_\pi(1) = \rho_\pi(2) = 0 \qquad{\rm and}\qquad \rho_\pi(3) = 1.$$
Furthermore, $\pi$ has a unique $G$-ascent (at~1),
and breaking $132$ after the~1
gives two subwords, one with two letters and the other with one letter,
so the $G$-ascent type of~$132$ is the integer partition~$(2,1)$.

\theorem\Gascent{Let $G$ be a graph with $V(G)=[n]$.  If
$$\X_G = \sum_\lambda N_\lambda \tilde\Xi_\lambda$$
is the $\tilde\Xi$-expansion of~$\X_G$,
then $N_\lambda$ is the number of permutations of~$[n]$ with $G$-ascent
type~$\lambda$.}

\proof
Given a stable partition~$\sigma$ of~$G$ and a linear ordering
$(\sigma_1, \sigma_2, \ldots, \sigma_\ell)$ of its blocks, we
will define an associated permutation of~$[n]$
(represented in one-line notation)
that will consist of a certain permutation of the elements of~$\sigma_1$
(the exact permutation will be specified in a moment) followed by
a certain permutation of the elements of~$\sigma_2$, and so on.

Before we specify how the elements within each~$\sigma_i$ are ordered,
we first notice that if $\pi_1$ and~$\pi_2$ are any two permutations
of the above form, then $\rho_{\pi_1}(v) = \rho_{\pi_2}(v)$ for all~$v$,
because there are no edges of~$G$ between vertices in the same block
of~$\sigma$.  Thus we may speak of the ``rank of a vertex'' without
ambiguity, even before specifying the permutation.

We now single out a particular permutation of the above form by arranging
the vertices within each block as follows: first arrange the vertices
in decreasing order of rank, and then within each rank, arrange the
vertices in numerically decreasing order.
The motivation for this arrangement
is that in the resulting permutation,
the only place that a $G$-ascent can occur
is between adjacent blocks of~$\sigma$.
This fact will play an important role shortly.

We next claim that the number of ordered stable partitions of type~$\mu$
that give rise to a fixed permutation~$\pi$ equals
$\ell(\mu)!\,a_{\lambda\mu}$, where $\lambda$ is the $G$-ascent type
of~$\pi$ and $a_{\lambda\mu}$ is the coefficient of $\tilde m_\mu$
in~$\tilde\Xi_\lambda$.  To see this, let $\sigma$ be
the ordered partition obtained
by breaking $\pi$ after each number at which there is a $G$-ascent.
A \itc{refinement} of~$\sigma$ is defined to be an ordered partition 
obtained by breaking $\pi$ at each break of~$\sigma$, plus (optionally)
any number of other locations.  From
the definition of~$\tilde\Xi_\lambda$,
we see that to prove our claim, it suffices to show
that it is precisely the type~$\mu$ refinements of~$\sigma$ that
generate~$\pi$.
We note first that it is only such refinements that can possibly
generate~$\pi$, by the remark at the end of the preceding paragraph.
It thus remains to show that every such refinement is in fact an
ordered \itc{stable} partition that generates~$\pi$.

To see that every refinement is stable, it suffices to show that
$\sigma$ is stable.  But if $i<j$ and $\pi(i)$ and~$\pi(j)$ are adjacent,
then the rank of $\pi(j)$ must exceed the rank of~$\pi(i)$, so there
must be at least one $G$-ascent between $i$ and~$j$, and thus $\pi(i)$
and $\pi(j)$ must be in different blocks of~$\sigma$.

To see that every refinement generates~$\pi$, note that the only way
this could fail to happen is if the elements in some block of the
refinement are not ordered properly (i.e., in decreasing order of
rank, and in decreasing numerical order within each rank).  However,
rank is defined without reference to partitions, so when we pass from
$\sigma$ to one of its refinements, the vertices must be ordered
properly in the refinement, provided that they were ordered properly
in~$\sigma$.  Finally, the vertices of~$\sigma$ are ordered properly,
because otherwise there would be a $G$-ascent inside a block of~$\sigma$.
This completes the proof of our claim.

Thus, if we let $b_\mu$ be the number of stable partitions of~$G$
of type~$\mu$, then for all $\mu$,
$$\ell(\mu)!\, b_\mu = \sum_\lambda N_\lambda \ell(\mu)!\, a_{\lambda\mu}.$$
Multiplying both sides by $\tilde m_\mu / \ell(\mu)!$ and summing
over~$\mu$,
$$\X_G = \sum_\mu b_\mu \tilde m_\mu =
  \sum_\lambda N_\lambda \sum_\mu a_{\lambda\mu} \tilde m_\mu
  = \sum_\lambda N_\lambda \tilde\Xi_\lambda.\qed$$

\sec The Poset Chain Conjecture revisited

In this section we give a new combinatorial proof of a (known) special
case of the Poset Chain Conjecture.  We hope that our method
can be generalized to handle more cases, although so far we have not
been able to do~so.

Let $P$ be a poset.
Recall from [GV2] or~[Ga1] that a \itc{$P$-array} is an array
$$\matrix{v_{1,1}&v_{1,2}&\ldots\cr
          v_{2,1}&v_{2,2}&\ldots\cr
          \ldots&&\cr}$$
where each $v_{i,j}$ is either undefined or an element of~$P$ and such
that

\item{1.} every element of~$P$ appears exactly once in the array, and

\item{2.} for all $i,j\ge1$, if $v_{i,j+1}$ is defined,
then $v_{i,j}$ is defined and
$v_{i,j} < v_{i,j+1}$ in~$P$.

The \itc{shape} of a $P$-array is the sequences of
the lengths of (the defined portion of) the rows.
A \itc{$P$-tableau} is a $P$-array such that

\item{3.} for all $i,j\ge1$, if $v_{i+1,j}$ is defined,
then $v_{i,j}$ is defined and
$v_{i+1,j} \not< v_{i,j}$ in~$P$.

We then have the following result of Gasharov~[Ga1].

\proposition\ptableaux{If $P$ is a $({\bf3+1})$-free poset and
$$\X_{G(P)} = \sum_\lambda a_\lambda s_\lambda$$
is the $s$-expansion of~$\X_{G(P)}$, then $a_\lambda$ is the
number of $P$-tableaux of shape~$\lambda$.\qed}

It is also well known [Mac, Table~1, p.~56] that
$$e_\lambda = \sum_\mu K_{\mu'\lambda} s_\mu,$$
where $K_{\mu'\lambda}$ is the number of semi-standard (i.e.,
row-nondecreasing and column-increasing) Young tableaux of shape~$\mu'$
and content~$\lambda$.  This points the way to a possible strategy
for proving $e$-positivity of a ($\bf 3+1$)-free poset~$P$ combinatorially:
try to find a partition of the set of all $P$-tableaux such that
for each block there is an integer partition~$\lambda$
such that the number of $P$-tableaux of shape~$\mu$ in that block
equals~$K_{\mu'\lambda}$.  From the above facts,
we see that the existence of such a partition
implies that $P$ is $e$-positive.

So far this method has not yielded new $e$-positivity results,
but it does lead to a combinatorial proof of the known fact~[S-S]
that $\bf 3$-free posets (i.e., posets which do not contain an
induced subposet isomorphic to a three-element chain) are $e$-positive.

\proposition\threefree{If $P$ is a $\bf 3$-free poset, then $\X_{G(P)}$
is $e$-positive.}

\proof
Let us define a \itc{$P$-diagram} to be an arrangement of the elements
of~$P$ into two columns, justified along the top edge, such that the
height of the right-hand column does not exceed the height of the
left-hand column.  (The columns are allowed to be empty.)
We define the \itc{popping} operation as follows.
To pop a $P$-diagram,
remove the bottommost element from the right-hand column
and place it at the bottom of the left-hand column.
(It is illegal to pop $P$-diagrams with empty right-hand columns.)

Assume now that $P$ is $\bf 3$-free.  Notice that this implies that
every $P$-array has at most two columns, so that in particular every
$P$-array is a $P$-diagram.  If $T$ and~$T'$ are $P$-tableaux, we
define the relation~$\sim$ by letting $T\sim T'$ if one of them
can be obtained from the other by a sequence of pops.  Clearly $\sim$
is an equivalence relation, and hence it induces a partition
$\pi = (\pi_1, \pi_2, \ldots)$ of the set of all $P$-tableaux.

For each~$i$, let $T_i$ be the $P$-tableau in~$\pi_i$
with the highest right-hand column.  Now fix any~$i$.
We claim that any $P$-diagram obtained from~$T_i$ by a sequence of pops
is in fact a $P$-tableau (and therefore lies in~$\pi_i$).  To see this,
note that the only way a problem can arise is if an element of~$P$ that
is shifted to the left-hand column in the course of a pop
turns out to be less than the element that it ends up sitting underneath.
However, because $T_i$ is a $P$-array, any element~$x$ in the right-hand
column of~$T_i$ is larger than at least one other element of~$P$ (namely,
the element immediately to the left of~$x$ in~$T_i$), and since $P$ is
$\bf 3$-avoiding, $x$ cannot be less than any other element in~$P$.
This establishes our claim.

Let $\lambda^i$ be the conjugate of the shape of~$T_i$.
In view of the remarks preceding the theorem statement, we see
that all there is left to prove is that for all~$i$, the number of
$P$-tableaux in~$\pi_i$ of shape~$\mu$ equals $K_{\mu'\lambda^i}$.
Now since $\lambda^i$ has at most two parts, $K_{\mu'\lambda^i} = 1$
if $|\mu'| = |\lambda^i|$,
$\mu'$ has at most two parts and $\mu'_1 \ge \lambda^i_1$, and
$K_{\mu'\lambda^i}=0$ otherwise.  Then from
the previous paragraph we see that there is indeed exactly
one $P$-tableau in~$\pi_i$ of shape~$\mu$ when $K_{\mu'\lambda^i}=1$
and there are no $P$-tableaux of other shapes.\qed

\sec Reconstruction

In [St2], several expansions of the chromatic polynomial are generalized
to expansions of~$\X_G$.  One notable omission is the multiplicative
expansion (see [Big, Chapter~11]), which has an important application
to the graph reconstruction conjecture.  It turns out, however, that
the multiplicative expansion of a graph invariant even more
general than~$\X_G$ was derived by Tutte a long time ago.
In this section we indicate the connection between Tutte's work
and the theory of the chromatic symmetric function, and in particular
we show how Tutte's work implies that $\X_G$ is reconstructible.

Recall that the \itc{list of vertex-deleted subgraphs} of a graph~$G$
is the multiset of unlabelled graphs
$$\{G\,\backslash v \mid v \in V(G)\}.$$
The graph~$G$ is \itc{reconstructible} if no other graph has the same
list of vertex-deleted subgraphs that $G$ does.  It is a major unsolved
problem in graph theory whether or not every graph with more than two
vertices is reconstructible.  See [Bon] for a recent survey.

If $G$ is a graph and $S\subset E(G)$, define $G\cdot S$ to be
the subgraph of~$G$ consisting of the edges of~$S$ together with
their incident vertices.
(In particular, $G\cdot S$ has no isolated vertices.)
Now for each isomorphism class~$K$
of connected graphs with at least two vertices,
let $x_K$ be an indeterminate.
Also let $t$ be an indeterminate distinct from all the~$x_K$.
Given any graph~$G$, define
$$J(G) \defeq \sum_{S\subset E(G)}
   \prod_{H\subset G\cdot S} t^{|V(H)|} x_H,$$
where the product is over all connected components~$H$ of~$G\cdot S$.
(If this product is empty we take its value to be one.)
We then have the following theorem of Tutte [Tu1,~6.6].

\proposition\tutte{The coefficients of all the terms of~$J(G)$ are
reconstructible, save possibly those terms containing $x_H$ where
$V(H)=V(G)$.\qed}

We remark, for the benefit of the reader who wishes to consult Tutte's
paper, that we have used $J(G)$ for what Tutte calls $J(E(G))$,
that we have used $x_H$ in place of~$f(H)$,
and that we are restricting our attention to the $C(G;1a)$~case.

\relax From \tutte\ it is easy to derive the following result.

\proposition\XGreconstruct{$\X_G$ is reconstructible.}

\proof
\relax From [St2, Theorem~2.5], we have
$$\X_G = \sum_{S\subset E(G)} (-1)^{|S|} p_{\lambda(S)},$$
where $\lambda(S)$ denotes the partition of~$d$ whose parts are equal
to the vertex sizes of the connected components of the spanning subgraph
of~$G$ with edge set~$S$.  It suffices to show that
$$\tilde\X_G \defeq \sum_{S\subset E(G)} (-1)^{|S|}
   \prod_{H\subset G\cdot S} p_{|V(H)|}$$
is reconstructible; what we have done is to ``set $p_1=1$ in~$\X_G$,''
which might appear at first to result in some loss of information,
but it is easy to see that $\X_G$ can be recovered from $\tilde\X_G$
using the fact that $\X_G$ is homogeneous.

Now set $t=1$ and $x_H = (-1)^{|E(H)|} p_{|V(H)|}$.  Then $J(G)$ becomes
$\tilde\X_G$, and from \tutte\ it follows that all the coefficients
in the $p$-expansion of~$\tilde\X_G$ can be reconstructed except possibly
for the coefficient of~$p_n$, where $n=|V(G)|$.  However, the coefficient
of~$p_n$ in~$\tilde\X_G$ equals the coefficient of~$p_n$ in~$\X_G$, which
in turn equals the coefficient of~$i$ in the chromatic
polynomial~$\chi_G(i)$ (as we can see by specializing
$\X_G(1^i)=\chi_G(i)$).  Since the chromatic polynomial is reconstructible
([Tu1, 7.5] or [Tu2]), this completes the proof.\qed

\sec Superfication

Recall that the \itc{superfication} of a symmetric function~$g$ is defined
by
$$g(x/y) \defeq \omega_y g(x,y).$$
Stanley~[St2] has briefly considered the superfication of~$\X_G$ and has
also asked what can be said about the two-variable polynomial
$\X_G(1^i\!/1^j)$.  We investigate the latter question in this section.

Some results can be obtained trivially by specializing known theorems
about~$\X_G$.  For example, we can specialize [St2, Theorem~4.3] as
follows.  Let $\Bbb P$ and $\bar{\Bbb P}$ be two disjoint copies of
the positive integers.  Denote the elements of~$\Bbb P$ by $1,2,3,\ldots$
and denote the elements of~$\bar{\Bbb P}$ by $\bar1,\bar2,\bar3,\ldots\,$.
Linearly order the disjoint union
${\Bbb P} \cup \bar{\Bbb P}$ by using the natural order on each of
the sets $\Bbb P$ and~$\bar{\Bbb P}$ and, additionally,
declaring every element of~$\Bbb P$ to be less than every element
of~$\bar{\Bbb P}$.  The following proposition then follows immediately
from [St2, Theorem~4.3].

\proposition\supercolor{For any graph~$G$, $\X_G(1^i\!/1^j)$ is
the number of pairs $({\Frak o},\kappa)$ such that $\Frak o$ is
an acyclic orientation of~$G$ and
$$\kappa:V(G)\to\{1,2,\ldots,i\}\cup \{\bar1,\bar2,\ldots,\bar j\}$$
is a map such that
(a)~if $u\to v$ is an edge of~$\Frak o$ then
$\kappa(u)\ge\kappa(v)$, and (b)~if $u\to v$ is an edge of~$\Frak o$
and both $\kappa(u)$ and~$\kappa(v)$ lie in~$\Bbb P$, then
$\kappa(u)>\kappa(v)$.\qed}

A second example is the following result, whose proof we omit since
it follows a standard line of argumentation in the the theory of
$P$-partitions and is somewhat long.  The reader is referred
to~[St3, Chapter~4] and~[St5] for the relevant theory and definitions.

\proposition\Ppartitions{For any graph~$G$ with $n$ vertices,
$$\X_G(1^i\!/1^j) = \sum_{\Frak o} \sum_{\pi\in{\Scr L}(\bar{\Frak o})}
  \sum_{k=0}^n {i + D_{n-k}(\pi) \choose n-k}
               {j + A_k(\pi) \choose k},$$
where the first sum is over all acyclic orientations~$\Frak o$ of~$G$,
$\bar{\Frak o}$ denotes the poset that is the transitive closure
of~$\Frak o$, ${\Scr L}$ denotes the Jordan-H\"older set, $D_k(\pi)$
denotes the number of descents in the first $k$ digits of~$\pi$, and
$A_k(\pi)$ denotes the number of ascents in the last $k$ digits
of~$\pi$.\qed}

For our third and final example we need a definition:
$$\tilde\chi_G(m,n) \defeq \X_G(1^{(m-n)/2}\!/1^{(m+n)/2}).$$

\proposition\mnpositive{The coefficients of the polynomial
$\tilde\chi_G(m,n)$ are nonnegative integers.}

\proof
Corollary~2.7 of~[St2] states that $\omega \X_G$ is $p$-positive, i.e.,
$\X_G$ is a nonnegative integer combination of~$(\sgn\lambda)p_\lambda$,
or a polynomial in the signed power sums $(\sgn d)p_d = (-1)^{d-1}p_d$
with nonnegative integer coefficients.  Now
$$(-1)^{d-1} p_d(x/y) = (-1)^{d-1}\bigl(p_d(x) + (-1)^{d-1}p_d(y)\bigr)
  = (-1)^{d-1} p_d(x) + p_d(y),$$
so
$$(-1)^{d-1} p_d(1^i\!/1^j) = (-1)^{d-1} i + j.$$
If we set $m=j+i$ and $n=j-i$, then we see that $\X_G(1^i\!/1^j)$ is
a polynomial in $m$ and~$n$ with nonnegative integer coefficients.
By rewriting the equations $m=j+i$ and $n=j-i$ in the form
$i=(m-n)/2$ and $j=(m+n)/2$ we see that the proposition follows.\qed

Apart from such specializations of known theorems,
however, it seems hard to find significant results.  For example,
it is easily checked that no deletion-contraction recurrence is possible.
The best we have been able to do is to find a recurrence for trees
which has some interesting properties.  Unfortunately we have not
found any applications for this recurrence, although it seems that
it might be possible to use it to obtain some partial results towards
the question, posed by Stanley, of whether nonisomorphic trees have
distinct chromatic symmetric functions.

To state the recurrence we must first define
two invariants that are closely related to
$\tilde\chi_G$.  Let $T$ be a \itc{rooted} tree, i.e., a tree with a
distinguished vertex.  Let $T'$ be the tree obtained by adjoining an
extra vertex~$v$ to~$T$ that is adjacent to the root of~$T$ (and to
no other vertices).  Define $\alpha_T(m,n)$ to be
the number of pairs $({\Frak o},\kappa)$
such that $\Frak o$ is an acyclic orientation of~$T'$ and
$\kappa:V(T')\to{\Bbb P}\cup\bar{\Bbb P}$ is a map that sends $v$ to~$1$,
whose image lies in
$$\{1,2,\ldots,(m-n)/2\} \cup \{\bar1,\bar2,\ldots,\overline{(m+n)/2}\},$$
and that also satisfies the conditions (a) and~(b) of \supercolor.
Define $\beta_T(m,n)$ in the same way except with the phrase
``sends $v$ to~$1$'' replaced by ``sends $v$ to~$\bar 1$.''
We then have the following theorem.

\theorem\treerecur{Let $T$ be a rooted tree with root~$v$, and let
$\Scr F$ be the family of rooted trees that results when $v$ is deleted
from~$T$.  (The root of a tree in~$\Scr F$ is the vertex that was
adjacent to~$v$ before deletion.)  Then
$$\eqalign{2\alpha_T(m,n) &=
  (m-n-2)\prod_{S\in\Scr F} \alpha_S(m,n) +
    (m+n)\prod_{S\in\Scr F} \beta_S(m,n)\cr
{\it and}\qquad 2\beta_T(m,n) &=
  (m+n+2)\prod_{S\in\Scr F} \beta_S(m,n) +
    (m-n)\prod_{S\in\Scr F} \alpha_S(m,n).\cr}$$}

\proof
Fix an arbitrary map $\kappa:V(T)\to\{1,2,\ldots,(m-n)/2\} \cup
\{\bar1,\bar2,\ldots,\overline{(m+n)/2}\}$.
We seek acyclic orientations~$\Frak o$ of~$T$ that are
compatible with $\kappa$ in the sense that $(\Frak o,\kappa)$
satisfies conditions (a) and~(b) of \supercolor.
Observe first that if $\kappa$ maps any two adjacent vertices
to the same element of~$\Bbb P$, then no acyclic orientations
are compatible with~$\kappa$.
Otherwise, we note that conditions (a) and~(b)
force a particular orientation of every edge except those
edges whose endvertices are mapped to the same element of~$\bar P$,
in which case the conditions (a) and~(b) impose no constraint on
the orientation.  Now since $T$ is a tree, \itc{every} orientation
of its edges is an acyclic orientation, and in particular every
orientation of the edges of~$T$ that meets the necessary conditions
just stated is in fact an acyclic orientation compatible with~$\kappa$.

With this in mind, we now consider $\alpha_T$.  Adjoin a vertex~$v$
to form~$T'$, and map $v$ to~$1$.  We wish to extend this to a
map~$\kappa$ of all the vertices of~$T'$ into
$$\{1,2,\ldots,(m-n)/2\} \cup \{\bar1,\bar2,\ldots,\overline{(m+n)/2}\},$$
and for each such $\kappa$ we want to find all compatible acyclic
orientations of~$T'$.  The total number of such compatible pairs will
give us $\alpha_T(m,n)$.  Begin by splitting into two cases:
in the first case, the root of~$T$ is mapped into~$\Bbb P$, and in
the second case, the root of~$T$ is mapped into~$\Bbb P$.  In the
first case, mapping the root to~1 results in adjacent vertices being
mapped to the same element of~$\Bbb P$, so we discard this possibility.
There remain
$${m-n\over 2} - 1 = {m - n - 2\over 2}$$
other possible colors in~$\Bbb P$ for the root of~$T$.  Now comes the
key observation: for each such coloring of the root of~$T$, we obtain
$$\prod_{S\in\Scr F} \alpha_S(m,n)$$
corresponding compatible pairs.  For once the color of the root
of~$T$ is fixed, each tree in~$\Scr F$ may be oriented and colored
independently, and the number of ways of doing this for a particular
$S\in\Scr F$ is just $\alpha_S(m,n)$---the root of~$T$ plays
the role of the adjoined vertex of~$S'$.  (There is a slight
technicality here in that in the definition of~$\alpha_S$
the adjoined vertex is required to be mapped to~$1$, whereas here
it may be mapped to an arbitrary element of~$\Bbb P$, but it is
clear that there is a bijection between the two sets of configurations,
with the only difference being that the orientation of the edge joining
the root of~$S$ to the root of~$T$ may need to be reversed in some cases.)
Furthermore, the orientation of the edge between $v$ and the root of~$T$
is forced.  Thus the total number of compatible pairs in the first case
is
$$\biggl({m-n-2\over2}\biggr) \prod_{S\in\Scr F} \alpha_S(m,n).$$

If we now consider the second case, we see that there are $(m+n)/2$ choices
of colors in~$\bar{\Bbb P}$ for the root of~$T$ and that each such choice
forces the orientation of the edge between $v$ and the root of~$T$.  The
number of ways of extending each such configuration to the entire tree is
$$\prod_{S\in\Scr F} \beta_S(m,n).$$
Adding up the two cases proves the first formula of our theorem.

The second formula is proved similarly; the only new twist is that
instead of \itc{discarding} the case where $v$ and the root of~$T$
are assigned the same color, we must \itc{count it twice,} since
in this case the edge between $v$ and the root of~$T$ may be oriented
in either direction.  This accounts for the term $(p+q+2)$.\qed

\treerecur\ allows us to state the precise relationship between
$\tilde\chi$ and $\alpha$ and~$\beta$.

\proposition\oneandtwo{Let $T$ be a rooted tree.  Then
$$\tilde\chi_T(m,n) = {(n+m)\alpha_T(m,n) + (n-m)\beta_T(m,n)
   \over 2(n+1)}.$$}

\proof
Let $\Scr F$ be the family of rooted trees resulting when the root of~$T$
is deleted.  Then it is clear from the definitions of $\tilde\chi$,
$\alpha$ and~$\beta$
$$\tilde\chi_T(m,n) = \biggl({m-n\over 2}\biggr)\prod_{S\in\Scr F}
   \alpha_S(m,n) + \biggl({m+n\over 2}\biggr)\prod_{S\in\Scr F}
   \beta_S(m,n),$$
where the two terms in this sum correspond to mapping the root of~$T$
into $\Bbb P$ and~$\bar{\Bbb P}$ respectively.  We can then invert
the formulas in \treerecur\ to solve for the products in terms of
$\alpha_T$ and $\beta_T$ (since the determinant of
$$\pmatrix{m-n-2&m-n\cr m+n&m+n+2\cr}$$
is not identically zero) to obtain the desired formula.\qed

Notice that the left-hand side of \oneandtwo\ is independent of the choice
of root, even though the terms on the right-hand side are not.

If we compute a few examples by hand, we quickly notice that
$\alpha$ and~$\beta$ are closely related to each other.
More precisely, we say that
two polynomials $\gamma_1(m,n)$ and $\gamma_2(m,n)$ form a
\itc{related pair} if the sign of the coefficient of any nonzero term
$m^r n^s$ in~$\gamma_1(m,n)$ equals $(-1)^{d-r}$ (where $d$ denotes
the degree of~$\gamma_1$) and if $\gamma_2$ can be obtained from
$\gamma_1$ by changing all the minus signs to plus signs.

\proposition\plusminus{Let $T$ be a rooted tree with $d$ vertices.
Then the degree of~$\alpha_T$ equals~$d$ and $\alpha_T$
and~$\beta_T$ form a related pair.}

\proof
That the degree is~$d$ follows easily (e.g., by induction).
For the second part of the proposition, we proceed by induction on~$d$.
Let $\Scr F$ be the family of rooted trees obtained by deleting
the root of~$T$.
Then for all $S\in\Scr F$,
$\alpha_S$ and $\beta_S$ form a related pair
by the induction hypothesis,
since each $S\in\Scr F$ has fewer than $d$ vertices.
It follows easily that
$$\prod_{S\in\Scr F} \alpha_S \qquad{\rm and}\qquad
  \prod_{S\in\Scr F} \beta_S$$
form a related pair, with degree $d-1$.  We may rewrite \treerecur\
in the form
$$\eqalign{\alpha_T &=
  \hbox{$1\over2$}m\biggl(\prod_{S\in\Scr F} \beta_S +
                \prod_{S\in\Scr F} \alpha_S\biggr) +
  \hbox{$1\over2$}n\biggl(\prod_{S\in\Scr F} \beta_S -
                \prod_{S\in\Scr F} \alpha_S\biggr) -
  \prod_{S\in\Scr F} \alpha_S\cr
\beta_T &=
  \hbox{$1\over2$}m\biggl(\prod_{S\in\Scr F} \beta_S +
                \prod_{S\in\Scr F} \alpha_S\biggr) +
  \hbox{$1\over2$}n\biggl(\prod_{S\in\Scr F} \beta_S -
                \prod_{S\in\Scr F} \alpha_S\biggr) +
  \prod_{S\in\Scr F} \beta_S\cr}$$
where we have omitted some of the $m$'s and $n$'s for brevity.

We now show that the terms in $\alpha_T$ have the appropriate sign.
Since $\prod\alpha_S$ and $\prod\beta_S$ form a related pair, their
sum contains only terms whose power of~$m$ differs from the highest
power of~$m$ (namely $d-1$) by an even integer, and moreover the
surviving terms are all nonnegative.  Hence the contribution
to $\alpha_T$ from the first summand has the correct signs.
For the second summand, note that $\prod\beta_S-\prod\alpha_S$
contains only terms that differ from $d-1$ by an \itc{odd}
integer, and that all these terms are nonnegative.  Multiplying
by~$n$ does not change any exponents of~$m$, so the contribution
from the second summand also has the correct signs.  Finally, it is
clear that $-\prod\alpha_S$ also contributes the correct signs.

To conclude the proof it suffices to show that $\beta_T$ can be
obtained by changing all the minus signs in $\alpha_T$
to plus signs.  As we
noted before, the first two summands consist entirely of
nonnegative terms, so the only thing we need to check is that
changing all minus signs to plus signs in $-\prod\alpha_S$ gives
$\prod\beta_S$, but this follows from the induction hypothesis.\qed

We conclude this section with the remark that
if we could show that $\alpha_T$ is always irreducible,
then it would follow that distinct rooted trees
always have distinct~$\alpha_T$'s.  For suppose we are given a
polynomial and are told that it equals $\alpha_T$ for some rooted
tree~$T$.  Then \plusminus\ lets us compute $\beta_T$, and then
inverting the formulas in \treerecur\ allows us to compute
$\prod\alpha_S$ and $\prod\beta_S$.  By irreducibility we can
then recover the $\alpha_S$, and by induction we may assume that
all the $S$'s may be reconstructed from the corresponding $\alpha_S$.
This would then allows us to reconstruct~$T$.
Unfortunately, proving irreducibility seems even harder than
proving the original conjecture!

\sec $\X_G(t)$

In [St4] Stanley considers briefly some generalizations of~$\X_G$,
including an invariant that he calls $\X_G(t)$.  Conceivably, an
entire thesis could be written about~$\X_G(t)$, but here we present
only the most basic facts, including one result that is mentioned
in~[St4] but whose proof is omitted there.

Let $G$ be a graph and let $\kappa$ be a map from $V(G)$ into the
positive integers.  Define
$$x^\kappa \defeq \prod_{v\in V(G)} x_{\kappa(v)}$$
and say that an edge is \itc{monochromatic} if $\kappa$ maps its
endvertices to the same integer.
Stanley~[St4] defines
$$\X_G(t) \defeq \sum_\kappa (1+t)^{m(\kappa)} x^\kappa,$$
where the sum is over all maps $\kappa$ from $V(G)$ into the positive
integers and $m(\kappa)$ denotes the number of monochromatic edges.
The motivation for this definition is that if we set $n$ of the $x$'s
equal to one and the rest equal to zero, the resulting two-variable
polynomial is equivalent to the coboundary polynomial (and therefore
the Tutte polynomial).  See [B-O] for more details.

\theorem\xgtexpansion{For any graph~$G$,
$$\X_G(t) = \sum_{S\subset E(G)} t^{|S|} p_{\lambda(S)}.$$}

\noindent
(Remark: here $\lambda(S)$ is as in [St2, Theorem~2.5] or as in the
proof of \XGreconstruct.)

\proof
Fix $k\ge0$ and consider the coefficient of~$t^k$ in~$\X_G(t)$.  The
maps~$\kappa$ that contribute to this coefficient are those 
with $m(\kappa)\ge k$.  Each such map contributes
$${m(\kappa)\choose k} x^\kappa.$$
Regard the binomial coefficient here as choosing $k$ of the $m(\kappa)$
monochromatic edges.  Then we see
that we can sum these contributions in another way: first list all
$k$-subsets $S\subset E(G)$, and then for each such $S$,
consider all maps~$\kappa$ that make every edge in~$S$
monochromatic.  A moment's thought reveals that
this sum is precisely $p_{\lambda(S)}$, and the theorem follows.\qed

\corollary\xgtforest{Let $G$ be a forest.  Then the coefficient of
each power sum in the power sum expansion of~$\X_G(t)$ is a monomial
in~$t$.}

\proof
If $S\subset E(G)$, let $G(S)$ be the subgraph of~$G$ with vertex
set~$V(G)$ and edge set~$S$.
Fix any integer partition~$\lambda$.
Since $G$ is a forest, $G(S)$ is also a forest for any subset~$S$,
and moreover the number of edges of~$S$
equals the number of vertices of~$G$
minus the number of connected components of~$G(S)$, i.e.,
$$|S| = |V(G)| - \ell\bigl(\lambda(S)\bigr).$$
This means that in the sum in \xgtexpansion, every term in~$p_\lambda$
has a coefficient of $t^{|V(G)|-\ell(\lambda)}$, independent of~$S$.
This proves the corollary.\qed

\theorem\xgtmobius{For any graph~$G$,
$$\X_G(t) = \sum_{\pi,\sigma\in L_G} (1+t)^{n(\pi)}
   \mu(\pi,\sigma) p_{\sigma},$$
where $L_G$ is the lattice of contractions of~$G$ and $n(\pi)$
is the number of edges whose endvertices lie in the same block of~$\pi$.}

\proof
For each $\sigma\in L_G$, define
$$\X_\sigma = \sum_\kappa x^\kappa,$$
where the sum is over all maps such that the monochromatic edges are
precisely those edges with both endvertices in the same block of~$\sigma$.
This is the same definition that Stanley makes in his proof of the
M\"obius function formula for~$\X_G$ [St2, Theorem~2.6]; as Stanley
observes, every map~$\kappa$ belongs to exactly one~$\X_\sigma$, so
$$p_\pi = \sum_{\sigma\ge\pi} \X_\sigma$$
for all $\pi\in L_G$.  By M\"obius inversion,
$$\X_\pi = \sum_{\sigma\ge\pi} \mu(\pi,\sigma) p_\sigma.$$
But clearly
$$\X_G(t) = \sum_{\pi\in L_G} (1+t)^{n(\pi)} \X_\pi,$$
and the theorem follows.\qed

\sec Counterexamples

We conclude by mentioning two counterexamples.  In [St2] it is shown that
$\X_G$ does not determine~$G$.  One might wonder what properties of~$G$
the chromatic symmetric function \itc{does} determine.  Does, for example,
$\X_G$ determine whether or not $G$ is planar?  The answer is no, because
the graphs

\vskip 0.3truein
\hskip 0.8truein
\psfig{figure=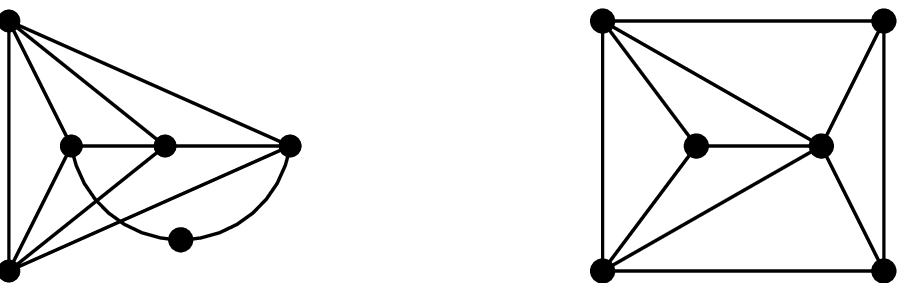}

\noindent have the same chromatic symmetric function
but the first graph is not planar (it is a
subdivision of~$K_5$) while the second one is.

Another question might be whether $\X_{G(P)}$ determines the dimension of
the poset~$P$.  (Recall that the dimension of a poset is the minimum number
of totally ordered sets needed to express the poset as an intersection of
totally ordered sets---see [Tro].)  Again, the answer is no; for example,
the incomparability graphs of the posets with Hasse diagrams

\vskip 0.3truein
\hskip 1truein
\psfig{figure=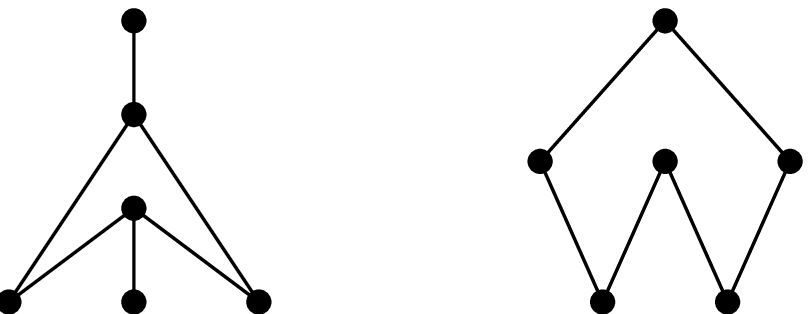}

\noindent
have the same chromatic symmetric function,
but the first poset has dimension~2 while the second poset has dimension~3.

\vfill\eject

\vbox to 1.5truein{\vfill}\leftline{\titlefont Bibliography}\vbox
        to 1.5truein{\vfill}

\parskip=12pt plus4pt minus4pt

\font\bigbold=cmbx12
\font\sc=cmcsc10
\def\sgn{{\rm sgn}\,}

\item{[Big]} {\sc N.~Biggs,} ``Algebraic Graph Theory,'' second edition,
Cambridge University Press, Cambridge, 1993.

\item{[Bon]} {\sc J.~A.~Bondy,} A graph reconstructor's manual,
\itc{in} ``Surveys in Combinatorics, 1991'' (Guildford, 1991),
London Math.\ Soc.\ Lecture Note Ser.~166,
Cambridge Univ.\ Press, Cambridge, 1991, pp.~221--252.

\item{[Bre]} {\sc F.~Brenti,} Expansions of chromatic polynomials and
log-concavity, \itc{Trans.\ Amer.\ Math.\ Soc.} {\bf 332} (1992), 729--756.

\item{[B-O]} {\sc T.~Brylawski and J.~Oxley,}
The Tutte polynomial and its applications, \itc{in}
``Matroid Applications'' (N.~White, Ed.),
Encyclopedia of Mathematics and Its Applications, Vol.~40,
Cambridge University Press, Cambridge, 1992, pp.~123--225.

\item{[BEGW]} {\sc J.~Buhler, D.~Eisenbud, R.~Graham, and C.~Wright,}
Juggling drops and descents, \itc{Amer.\ Math.\ Monthly}
{\bf 101} (1994), 507--519.

\item{[B-G]} {\sc J.~Buhler and R.~Graham,} A note on the binomial
drop polynomial of a poset, \itc{J.~Combin.\ Theory (A)} {\bf 66}
(1994), 321--326.

\item{[CSV]} {\sc L.~Carlitz, R.~Scoville, and T.~Vaughan,}
Enumeration of pairs of
sequences by rises, falls and levels, \itc{Manuscripta Math.} {\bf 19}
(1976), 211--243.

\item{[Cho]} {\sc T.~Chow,} The path-cycle symmetric function of a digraph, 
\itc{Advances in Math.,} in press.

\item{[CG1]} {\sc F.~Chung and R.~Graham,} On digraph polynomials, preprint
dated December~3, 1993.

\item{[CG2]} {\sc F.~Chung and R.~Graham,} On the cover polynomial of a
digraph, \itc{J.~Combin.\ Theory (B),} in press.

\item{[D-K]} {\sc F.~N.~David and M.~G.~Kendall,}
Tables of symmetric functions.~I,
\itc{Biometr.} {\bf 36} (1949), 431--449.

\item{[Dou]} {\sc P.~Doubilet,} On the foundations of combinatorial theory.~VII:
Symmetric functions through the theory of distribution and occupancy, 
\itc{Studies in Applied Math.} {\bf 51} (1972), 377--396.

\item{[Dwo]} {\sc M.~Dworkin,} personal communication.

\item{[GV1]} {\sc E.~Gansner and K.-P.~Vo,} The chromatic generating function,
\itc{Lin.\ and Multilin.\ Alg.} {\bf 22} (1987), 87--93.

\item{[G-J]} {\sc M.~R.~Garey and D.~S.~Johnson,} ``Computers and
Intractability: A Guide to the Theory of NP-Completeness,'' W.~H.~Freeman
and~Co., San Francisco, CA, 1979.

\item{[Ga1]} {\sc V.~Gasharov,} Incomparability graphs of $(\bf 3+1)$-free
posets are $s$-positive, \itc{in} Abstracts of the Sixth Conference on
Formal Power Series and Combinatorics, May 23--27, 1994.

\item{[Ga2]} {\sc V.~Gasharov,} On Stanley's chromatic symmetric function
and stable set polynomials, preprint dated January~8, 1995.

\item{[Ge1]} {\sc I.~M.~Gessel,} personal communication.

\item{[Ge2]} {\sc I.~M.~Gessel,}
Multipartite $P$-partitions and inner products of
skew Schur functions, \itc{in} ``Combinatorics and Algebra'' (C.~Greene,
Ed.), Contemporary Mathematics Series, Vol.~34, Amer.\ Math.\ Soc.,
Providence, R.I., 1984, pp.~289--301.

\item{[GV2]} {\sc I.~M.~Gessel and G.~X.~Viennot,}
Determinants, paths, and plane partitions, preprint dated July 28, 1989.

\item{[Gol]} {\sc J.~R.~Goldman,} personal communication.

\item{[GJW1]} {\sc J.~R.~Goldman, J.~T.~Joichi, and D.~E.~White,}
Rook theory~I. Rook equivalence of Ferrers boards,
\itc{Proc.\ Amer.\ Math.\ Soc.} {\bf 52} (1975), 485--492.

\item{[GJRW]} {\sc J.~R.~Goldman, J.~T.~Joichi, D.~L.~Reiner,
and D.~E.~White,} Rook theory.~II: Boards of binomial type,
\itc{SIAM J.~Applied Math.} {\bf 31} (1976), 618--633.

\item{[GJW2]} {\sc J.~R.~Goldman, J.~T.~Joichi, and D.~E.~White,}
Rook theory~III. Rook polynomials and the chromatic structure of graphs,
\itc{J.~Combin.\ Theory (B)} {\bf 25} (1978), 135--142.

\item{[GJW3]} {\sc J.~R.~Goldman, J.~T.~Joichi, and D.~E.~White,}
Rook theory---IV. Orthogonal sequences of rook polynomials,
\itc{Studies in Applied Math.} {\bf 56} (1977), 267--272.

\item{[GJW4]} {\sc J.~R.~Goldman, J.~T.~Joichi, and D.~E.~White,}
Rook theory.~V. Rook polynomials, M\"obius inversion and the umbral calculus,
\itc{J.~Combin.\ Theory (A)} {\bf 21} (1976), 230--239.

\item{[Gre]} {\sc C.~Greene,} Proof of a conjecture on immanants of
the Jacobi-Trudi matrix, \itc{Lin.\ Alg.\ Appl.} {\bf 171} (1992), 65--79.

\item{[Kre]} {\sc G.~Kreweras,} The number of more or less ``regular''
permutations, \itc{Fibonacci Quart.} {\bf 18} (1980), 226--229.

\item{[Lin]} {\sc N.~Linial,} Graph coloring and monotone functions on posets,
\itc{Disc.\ Math.} {\bf 58} (1986), 97--98.

\item{[Mac]} {\sc I.~G.~Macdonald,} ``Symmetric Functions and Hall
Polynomials,'' Oxford University Press, Oxford, 1979.

\item{[Rio]} {\sc J.~Riordan,} ``An Introduction to Combinatorial Analysis,''
John Wiley \& Sons, New York, 1958.

\item{[R-R]} {\sc S.~M.~Roman and G.-C.~Rota,} The umbral calculus,
\itc{Advances in Math.} {\bf 27} (1978), 95--188.

\item{[R-P]} {\sc M.~Rumney and E.~J.~F.~Primrose,} A sequence connected
with the sub-factorial sequence (Mathematical notes 3207), \itc{Math.\
Gazette} {\bf 52} (1968), 381--382.

\item{[Slo]} {\sc N.~J.~A.~Sloane,} ``The On-Line Encyclopedia of Integer
Sequences.''  [For information on accessing this encyclopedia, see
{\sc N.~J.~A.~Sloane,} An on-line version of the encyclopedia of integer
sequences, \itc{Electronic J.~Combin.} {\bf 1} (1994), F1.]

\item{[St1]} {\sc R.~P.~Stanley,} personal communication.

\item{[St2]} {\sc R.~P.~Stanley,} A symmetric function generalization of the
chromatic polynomial of a graph, \itc{Advances in Math.} {\bf 111} (1995),
166--194.

\item{[St3]} {\sc R.~P.~Stanley,}
``Enumerative Combinatorics,'' vol.~1, Wadsworth
\& Brooks/{\hskip0.01pt}Cole, Pacific Grove, CA, 1986.

\item{[St4]} {\sc R.~P.~Stanley,} Graph colorings and related symmetric
functions: ideas and applications, preprint.

\item{[St5]} {\sc R.~P.~Stanley,} Ordered structures and partitions, 
\itc{Memoirs Amer.\ Math.\ Soc.} {\bf 119} (1972).

\item{[S-S]} {\sc R.~P.~Stanley and J.~Stembridge,}
On immanants of Jacobi-Trudi
matrices and permutations with restricted position, \itc{J.~Combin.\
Theory (A)} {\bf 62} (1993), 261--279.

\item{[Tro]} {\sc W.~T.~Trotter,} ``Combinatorics and Partially Ordered
Sets: Dimension Theory,'' Johns Hopkins University Press, Baltimore, MD,
1992.

\item{[Tu1]} {\sc W.~T.~Tutte,} All the king's horses (a guide to
reconstruction), \itc{in} ``Graph Theory and Related Topics''
(J.~A.~Bondy and U.~S.~R. Murty, Eds.),
Proceedings of the conference held in honor of Professor W.~T.~Tutte
on the occasion of his sixtieth birthday, Waterloo, July 5--9, 1977, 
Academic Press, New York, NY, 1979, pp.~15--33.

\item{[Tu2]} {\sc W.~T.~Tutte,} The reconstruction problem in graph theory,
\itc{Brit.\ Polymer~J.,} September 1977, 180--183.

\item{[Vo]} {\sc K.-P.~Vo,}
Graph colorings and acyclic orientations, \itc{Lin.\
and Multilin.\ Alg.} {\bf 22} (1987), 161--170.

\bye